\documentclass[11pt,a4paper]{scrartcl} 
\pdfoutput=1
\usepackage[utf8]{inputenc} %input encoding
\usepackage{amstext} 
\usepackage{amsmath} % added by UM
\usepackage{amssymb} % added by UM
\usepackage[pdftex,plainpages=false]{hyperref}
\usepackage{verbatim} %added by UM
\usepackage[pdftex]{graphicx}

\newcommand{\eq}[1]{(\ref{eq:#1})}   
\newcommand{\Lb}[1]{\label{eq:#1}}   
\usepackage{subfig} 
\renewcommand{\d}{\,\mathrm{d}} 
\newcommand{\e}{\,\mathrm{e}} 
\renewcommand{\i}{\,\mathrm{i}\,}

\newcommand{\half}{\frac{1}{2}}

\newcommand{\N}{\mathbb{N}}

\newcommand{\R}{\mathbb{R}}     
\newcommand{\C}{\mathbb{C}}

\newcommand{\Ts}{\mathcal{T}}

\newcommand{\Vs}{\mathcal{H}}

\newcommand{\ep}{\; .}  

\newcommand{\ec}{\; , \quad} 
\newcommand{\elli}{,\ldots ,}

\renewcommand{\u}[1]{\underline{#1}} 

\newcommand{\timedot}{\centerdot}
\renewcommand{\dot}[1]{\overset{\timedot}{#1}}
\renewcommand{\ddot}[1]{\overset{\timedot\timedot}{#1}}

\newcommand{\id}{\mathbf{1}}

\newcommand{\norm}[1]{\parallel \! #1 \! \parallel}
\newcommand{\map}[2]{#1 \: \rightarrow \: #2 }
\newcommand{\set}[2]{\{#1 \; : \; #2\}}  
\newcommand{\mapto}[2]{#1~\,~\mapsto~\,~#2 }
\newcommand{\namedmapto}[3]{#1 \: : \: #2 \: \mapsto \: #3 }

\newcommand{\fullmapeq}[6]
{\begin{equation}\Lb{#1}
\begin{split}
#2 \: : \: #3 \: &\rightarrow \: #4 \\
#5 \: &\mapsto \: #6 \ep
\end{split}
\end{equation}}
\newcommand{\fullmapeqnp}[6]
{\begin{equation}\Lb{#1}
\begin{split}
#2 \: : \: #3 \: &\rightarrow \: #4 \\
#5 \: &\mapsto \: #6
\end{split}
\end{equation}}

\newcommand{\dglx}[2]{F( #2,\,#1 )}
\newcommand{\dgl}[2]{\dglx{ #1(#2)}{#2}}
\newcommand{\dglrhs}[1]{\dgl{\psi}{#1}}
\newcommand{\dgli}[1]{\dglx{\psi_{#1}}{t_{#1}}}
\newcommand{\dglrhst}{\dglrhs{t}}
\newcommand{\LF}{\mathcal{L}}
\newcommand{\ALF}{\mathcal{A}}
\newcommand{\stage}[1]{(t_{#1},\,\psi_{#1})}

\begin{document} 

\date{}     
\title{\vspace{0 mm} An asynchronous leapfrog method II}
\author{Ulrich Mutze
\footnote{www.ulrichmutze.de}}\label{fnAuthor}
\maketitle 
\begin{abstract}
A second order explicit one-step numerical method for the initial value 
problem of the general
ordinary differential equation is proposed. 
It is obtained by natural modifications of the well-known
leapfrog method, which is a second order, two-step, explicit method. 
According to the latter method, the input data for
an integration step are two system states, wich refer to different 
times. 
The usage of two states instead
of a single one can be seen as the reason for the robustness of the 
method. 
Since the time step size thus is part of the 
step input data, it is complicated to change this size 
during the computation of a discrete trajectory.
This is a serious drawback when one needs to implement automatic
time step control.

The proposed modification transforms one of the two input states into a 
velocity
and thus gets rid of the time step dependency in the step input data.
For these new step input data, the leapfrog method gives a unique
prescription how to evolve them stepwise.

The stability properties of this modified method are the same as 
for the original one: the set of absolute stability is the interval 
$[-\i,+\i]$ on the imaginary axis. 
This implies exponential growth of trajectories in situations where the 
exact trajectory has an asymptote.

By considering new evolution steps that are composed of two consecutive
old evolution steps we can average over the velocities of the sub-steps 
and get an integrator with a much larger set of absolute stability, which 
is immune to the asymptote problem. 

The method is exemplified with the equation of motion of a one-dimensional
non-linear oscillator describing the radial motion in the Kepler problem.
For this system the exact solution is cheaply computable.
Depending on the selected initial conditions, widely varying degrees of 
stiffness hold. 
For this system the proposed leapfrog methods are compared with 
the original leapfrog method, three explicit second order Runge-Kutta 
methods, and a St\o rmer-Verlet method.
In all tested cases the asynchronous leapfrog method was more accurate 
than the original method, in some cases considerably so.

A `numerical interaction picture' is introduced, in which the evolution of 
a dynamical system and its evolution as brought about by a numerical 
integrator are treated in analogy to the dynamics with and without 
interaction in the interaction picture of quantum mechanics.
The trajectories of a system in this representation are indicative for the 
accuracy of an integrator and provide a kind of a fingerprint of 
it.

For non-stiff initial conditions, the asynchronous leapfrog method turns 
out to be the most accurate among the methods under consideration. 
The averaged version behaves very robust for stiff initial conditions
since it is not reversible and slightly dissipative.
It is well suited as a general-purpose method for rough real-world 
problems.

Finally a automatic step control scheme is demonstrated. 
It works for both the asynchronous leapfrog methods and the Runge-Kutta 
methods. Here, averaging in the leapfrog method makes no difference and 
leapfrog methods are by a factor of about 4 more accurate than the 
Runge-Kutta methods.   
\end{abstract}

AMS Subject Classification (MSC2010): 65L04, 65L05, 65L07, 65L12, 
65P10, 65Q10


\section{Preface}
\label{preface}
This work originated from an expanded form of \cite{leapfrog a} which is 
accessible from my homepage as \cite{leapfrog b}. 
A forum discussion in \cite{rg} based on \cite{leapfrog b} stressed the 
instability properties of leapfrog methods and convinced me that the 
asynchronous leapfrog integrator is unsuitable as a general-purpose 
workhorse method unless its stability can be improved.
Most instructive in this respect were the Lotka-Volterra equations 
proposed to me by H.E. Lehtihet.  
Fortunately the densified form of the asynchronous leapfrog integrator 
suggests an computationally cheap averaging step which improves stability 
drastically. Definition of this new integrator and demonstration of its 
properties is the main topic of this work. 
Although meanwhile I worked out many applications to higher-dimensional 
problems and, with the spherical pendulum and the spinning top, also 
two with non-trivial geometric structure of its configuration space, I 
decided to include no new systems and to make only minor changes to
sections 3 to 5 of \cite{leapfrog b}. 
The material concerning the new integrator and a study of stability in all 
methods under consideration is added as two following sections. 
Of course, Abstract, Introduction, Acknowledgments, 
and References are changed to fit the present form of the article.

\section{Introduction}
\label{intro}
We consider the initial value problem of the general ordinary differential 
equation
\begin{equation}
\Lb{ev}
\dot{\psi}(t) = \dglrhst
\end{equation}
for a time-dependent quantity $\psi$ which takes values in a real 
finite-dimensional
vector space $\Vs$.
Here $F$ is a function $\map{\R \times \Vs}{\Vs}$ \footnote{
As is well-known, one may increase formal simplicity by transforming away 
the explicit 
$t$-dependence of the right-hand side of this equation, thus rendering it 
\emph{autonomous}. 
However, I refrain from assuming autonomy, since the algorithms to be 
considered should apply to 
time-dependent real-world problems directly, without a need to transform 
them into autonomous
form\@.}~\@.
Equations of this kind arise in various contexts:
\begin{enumerate}
 \item Ordinary differential equations of finite order for real or complex 
variables.
 \item Dynamical systems on finite-dimensional real or complex manifolds; 
in this case $\Vs$ is
 a space of coordinate values (the co-domain of the charts from a suitable 
atlas), and $F$ is the
coordinate representation of a vector field. 
 \item Discrete approximations to partial differential equations like the 
time-dependent
Schr\"{o}dinger equation or Maxwell's equations.
\end{enumerate}
Although the parameter $t$ will always be referred to as `time', it could 
have a different
meaning as in the following examples: length of a rod as a function of 
temperature, position
of a point on a curve as a function of the arc length separating it from 
some reference point
on the curve.

The \emph{computational initial value problem} associated with this 
equation \eq{ev} asks for an
algorithm which determines for each $\R$-valued increasing list $t_0, t_1 
\elli t_n$
and each value $\psi_0 \in \Vs$ (the initial value) a $\Vs$-valued list 
$\psi_1 \elli \psi_n$ such 
that the $\R \times \Vs$-valued list $(t_0,\psi_0), (t_1, \psi_1) \elli 
(t_n, \psi_n)$
is a reasonable approximation to a solution curve $\mapto{t}{\psi(t)}$, 
$\psi(t_0) = \psi_0$, 
of the differential equation
\eq{ev} whenever the regularity properties of $F$ suffice
for determining such a curve, and the gaps between adjacent $t-$values are 
small enough.
If such an algorithm works only for equidistant time lists 
( for which $t_i - t_{i-1}$ is independent of $i$ by definition) it is
said 
to be \emph{synchronous}
and otherwise it is said to be \emph{asynchronous}.
Asynchronous algorithms may be developed into adaptive ones, which adjust 
their step size
$t_{i+1}-t_{i}$ to the size (according to a suitable notion of size in 
$\Vs$)
of $\dglx{\psi_i}{t_i}$.
For a $\R$-valued function $F$ that depends on its second argument 
trivially, the initial value 
problem is simply the problem of computing the definite integral. 
It is straightforward and instructive to specialize the proposed
algorithms 
to this simplified concrete situation.

Starting from the well-known leapfrog algorithm, the present article 
develops and analyzes an economic and robust
asynchronous solution of the computational initial value problem
associated 
with \eq{ev}.
Section~\ref{leapfrog} recalls the leapfrog method, and 
Section~\ref{asynchronous} carries out the general development
of the new asynchronous algorithm.
Section~\ref{Kepler}  will apply this integration method and related 
methods to the differential equation defined by \eq{radial}.
Section~\ref{adalf} defines the averaged form of the asynchronous leapfrog 
integrator and analyzes the stability of it and of the 
second order explicit Runge-Kutta methods.
Section~\ref{Kepler 2} tests the new methods by applying them to the 
differential equation defined by \eq{radial}. 

\section{The leapfrog method}
\label{leapfrog}
\begin{figure}
\centering
\mbox{
\subfloat[Dependence of the accuracy on the initialization method.]
{\includegraphics[width=70mm]{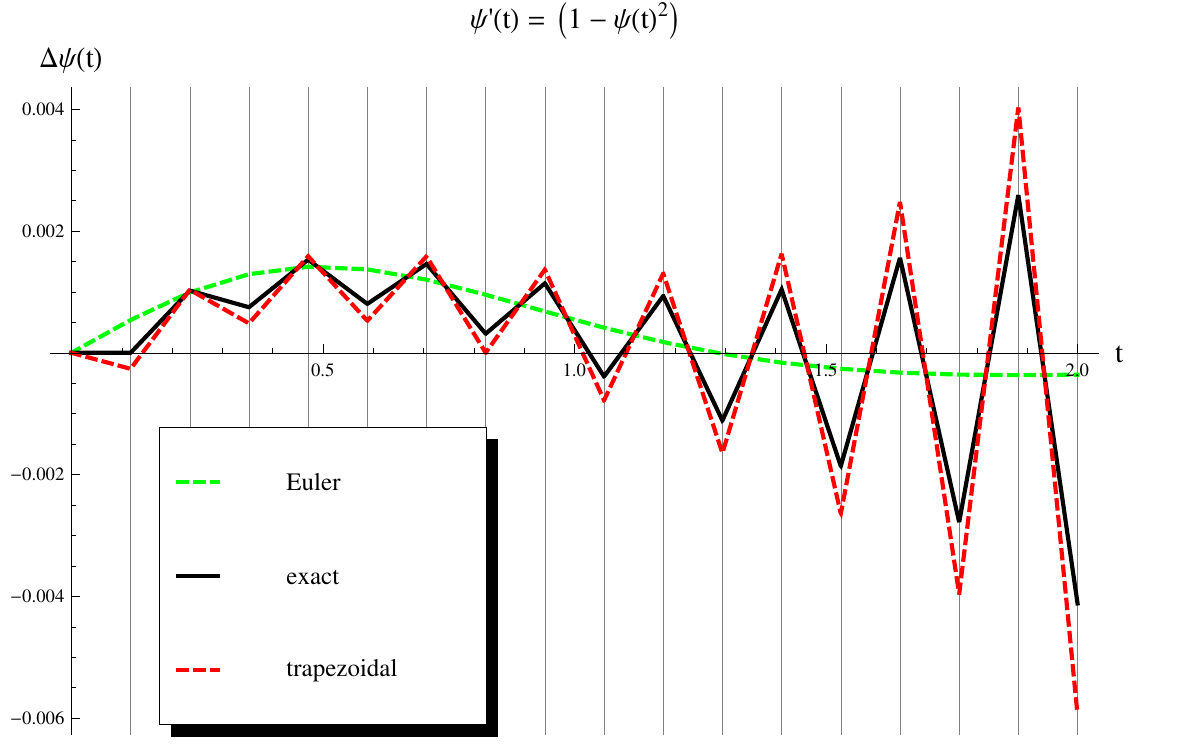}}
\quad
\subfloat[Exact (green) and discrete trajectory (red, Euler 
initialization).]
{\includegraphics[width=70mm]{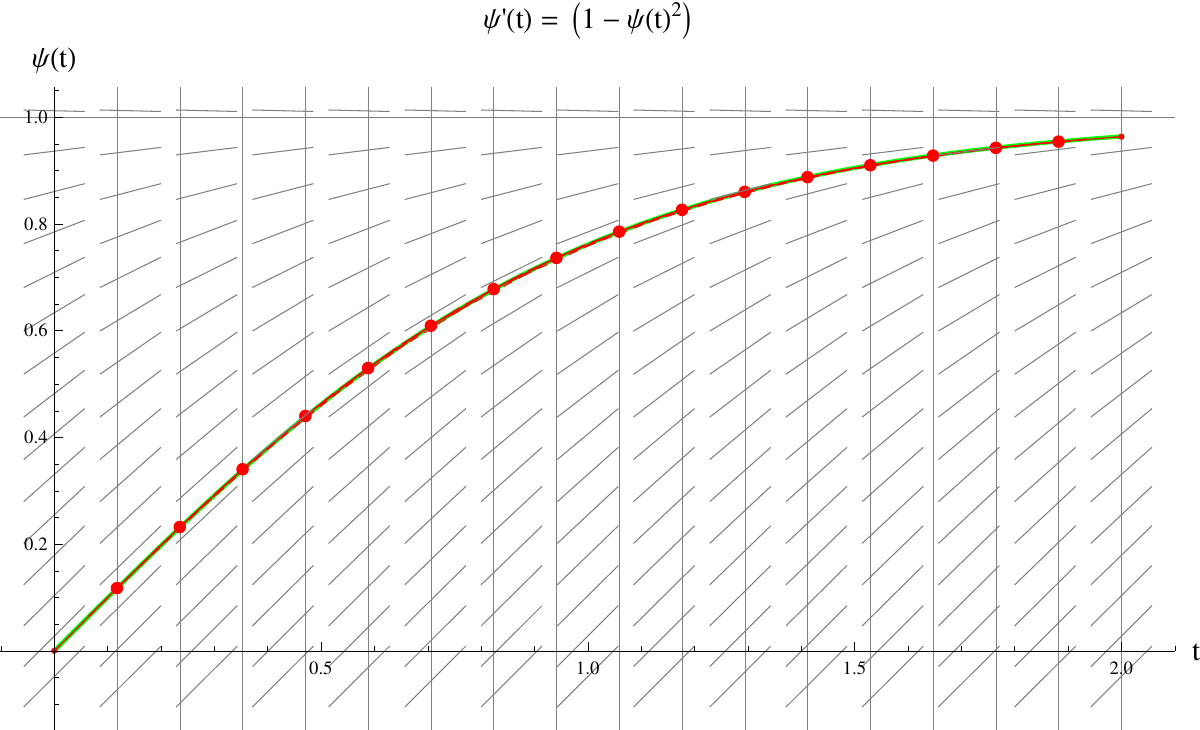}}
}
\caption{Case in which Euler initialization is superior.}
\label{fig:Euler good}
\end{figure}

A marvelously simple synchronous solution algorithm for the computational 
initial value problem 
of \eq{ev} is the \emph{leapfrog method}
or \emph{explicit midpoint rule}, see e.g.  \cite{Stuart}, eq. (3.3.11). 
It seems to be the first method that has been successfully applied to the
initial value problem of the time-dependent Schr\"{o}dinger equation (see 
\cite{askar} and the
citation of this work in \cite{talezer}).
It is most conveniently considered as the map
\fullmapeqnp{leapfrog0}{\LF}{(\R \times \Vs ) \times (\R \times \Vs )}{(\R 
\times \Vs ) 
\times (\R \times \Vs )}{((t_0,\psi_0), 
(t_1,\psi_1))}{((t_1,\psi_1),(t_2,\psi_2)) \ec}
where
\begin{equation}
\Lb{leapfrog}
\begin{split}
t_2 &:= 2 t_1 - t_0 \ec \\
\psi_2 &:= \psi_0 + (t_{2}-t_{0}) \,\dglx{\psi_1}{t_1}\ep
\end{split}
\end{equation}
The equivalent form 
\begin{equation}
\Lb{leapfrog2}
\begin{split}
\frac{t_0 + t_2}{2} &= t_1 \ec \\
\frac{\psi_2 - \psi_0}{t_2-t_0} &= \dglx{\psi_1}{t_1}
\end{split}
\end{equation}
of these equations, together with equation \eq{ev}, makes the reasons for 
choosing them evident:
\begin{equation}
\Lb{derivation}
\dglrhs{t_1} = \dot{\psi}(t_1) = \frac{\psi(t_2) - \psi(t_0)}{t_2-t_0} + 
O((t_2-t_0)^2) \ep
\end{equation}
Iterating the map $\LF$ determines a \emph{leapfrog trajectory}, 
$((t_j,\psi_j))_{j \in \N}$ 
if $(t_0,\psi_0)$ and $(t_1,\psi_1)$
are given:
\begin{equation}
\Lb{iteration}
(t_{i+1},\psi_{i+1}) = \pi_2( \LF ((t_{i-1},\psi_{i-1}),
(t_{i},\psi_{i}))) 
= 
\pi_2( {\LF}^{i} ((t_{0},\psi_{0}), (t_{1},\psi_{1}))) \ec
\end{equation}
where $\pi_2$ is the canonical projection to the second component of a
pair. 
According to the initial value problem we are given $t_0,\, t_1$ (which 
determines,
due to the assumed synchrony, all further $t$-values) and $\psi_0$. 
For starting iteration \eq{iteration} we need also $\psi_1$. 
This value has to be set in a way consistent with \eq{ev}, e. g. by
employing
the \emph{explicit Euler rule}
\begin{equation}
\Lb{euler}
\psi_1 := \psi_0 + (t_1-t_0) \,\dglx{\psi_0}{t_0} \ep
\end{equation}
One could expect that being more accurate here
would improve the accuracy of the leapfrog trajectory.
The optimum definition could be expected to be
\begin{equation}
\Lb{exact}
\psi_1 := \psi(t_1) \ec
\end{equation}
where $\psi(\cdot)$ is the \emph{exact} trajectory determined by \eq{ev}
and $\psi(t_0) = \psi_0$.
In a reasonably posed problem $t_1$ is close to $t_0$ and the exact 
trajectory up to $t_1$ can be arbitrarily well approximated by
sufficiently many steps of any consistent numerical method.
An intermediate position between methods \eq{euler} and \eq{exact} is
given by employing  the implicit \emph{trapezoidal rule}
\begin{equation}
\Lb{initialization implicit}
\psi_{1} = \psi_0 + (t_1-t_0) \, 
\frac{\dglx{\psi_0}{t_0}+\dglx{\psi_{1}}{t_1}}{2} \ec
\end{equation}
which can be efficiently solved by the iteration
\begin{equation}
\Lb{initialization iteration}
\psi^{(i+1)}_{1} := \psi_0 + (t_1-t_0) \, 
\frac{\dglx{\psi_0}{t_0}+\dglx{\psi^{(i)}_{1}}{t_1}}{2} \ec \psi^{(0)}_{1} 
:= \psi_0 \ep
\end{equation}
Let us study the behavior of these methods for a well-known non-linear 
differential equation for which the exact solution is known:
\begin{equation}
\Lb{tanh}
\dot{\psi}(t) = 1 - \psi(t)^2 \ec \psi(0) = 0 \quad \text{thus} \quad 
\psi(t) = \tanh(t) \ep
\end{equation}
Figure~\ref{fig:Euler good} shows the difference between the exact 
trajectory and the leapfrog
trajectories according to the three initialization methods introduced 
above.
Surprisingly the less sophisticated method \eq{euler} works best in this 
case.
% Very convenient structure for information-rich illustration! Needs 
%subfigure package
Unfortunately this is not the case in all interesting applications.
An example for this is shown in Figure~\ref{fig:Euler bad} for a
different differential equation
\footnote{The 'exact' solution is here the one provided by the default 
integrator of Mathematica.}
, where method  \eq{euler} gives by far the largest error. 

In a sense, selecting $\psi_1$ is a degree of freedom available for
solving the initial value problem and extending it to `hopping solutions':
\begin{figure}
\centering
\mbox{
\subfloat[Dependence of the accuracy on the initialization method.]
{\includegraphics[width=72mm]{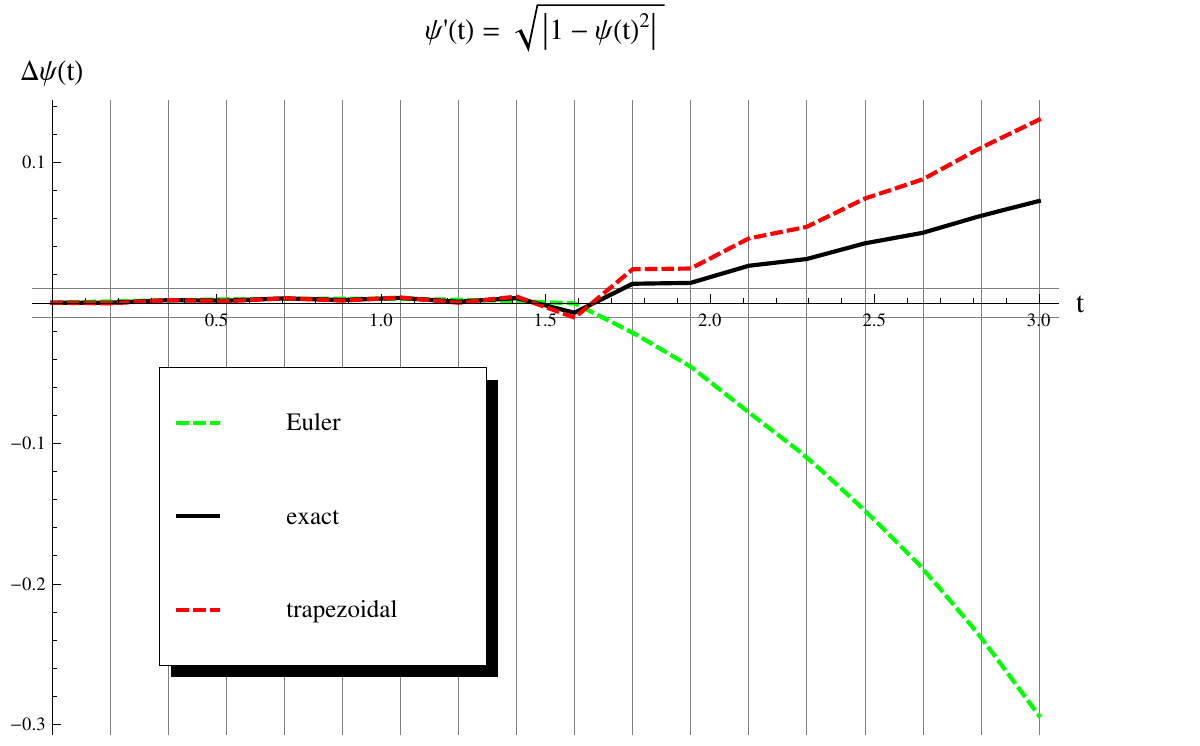}}
\quad
\subfloat[Exact (green) and discrete trajectory (red, Euler 
initialization).]
{\includegraphics[width=72mm]{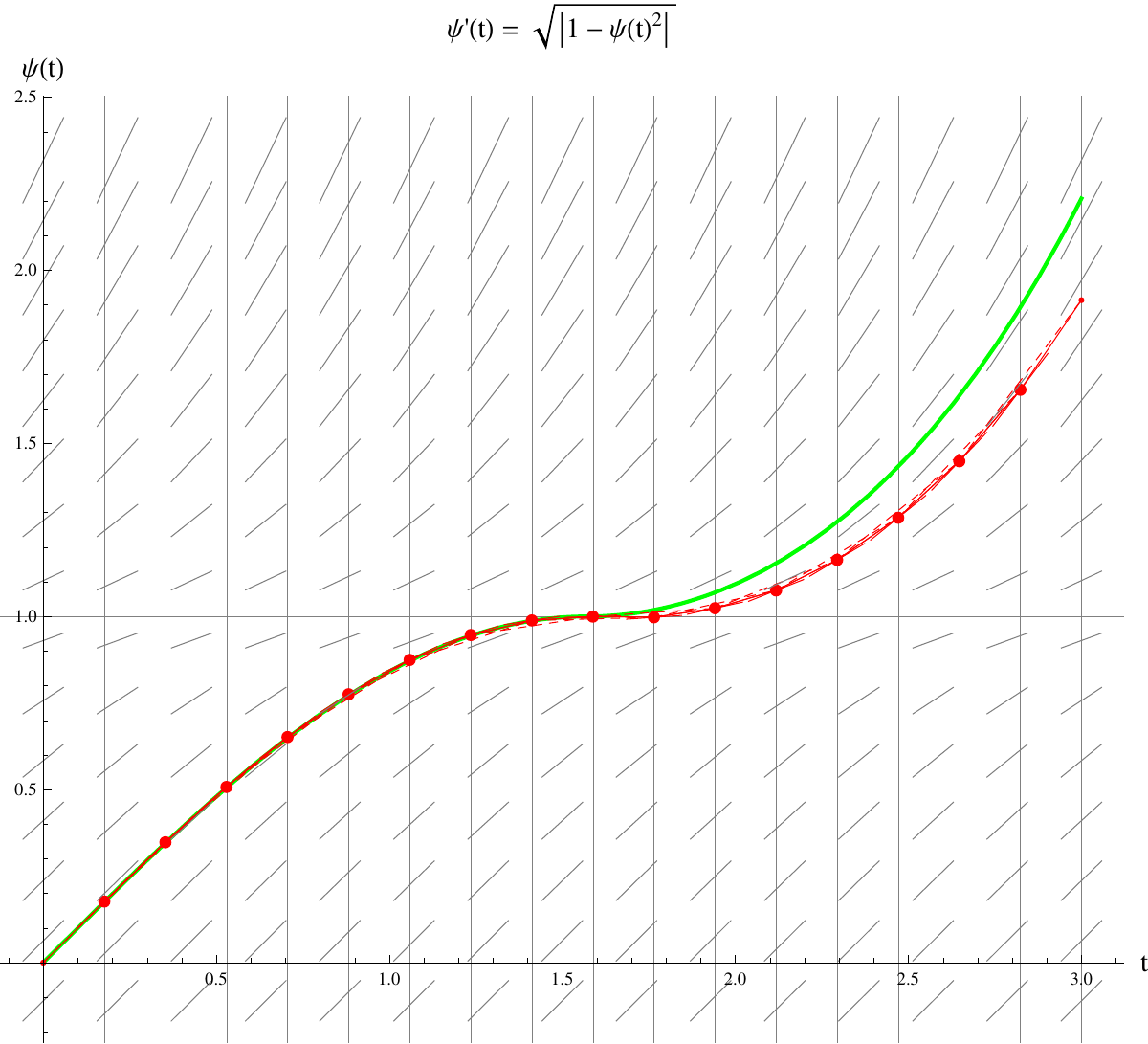}}
}
\caption{Case in which Euler initialization is inferior.}
\label{fig:Euler bad}
\end{figure}
%Figure~\ref{fig:Euler good}(b) is interesting.
Even for  arbitrary $\psi_1$, the iteration \eq{iteration} can be 
used to define a leapfrog trajectory.
If $\psi_1$ is not too far away from the value determined by
\eq{euler}, the leapfrog trajectory is a zigzag line which
tends to wiggle around a trajectory of differential equation \eq{ev} and
initial data $(t_0,\psi_0)$. 
See Figure~\ref{fig:zigzag}, where the leapfrog trajectory of 
Figure~\ref{fig:Euler good}(b)
is modified by shifting $\psi_1$ considerably from the position
determined by \eq{euler}.
It suggests a growing uncertainty in position and shows some similarity
with broadening wave packets.
The map \eq{leapfrog0}, considered as a discrete dynamical system, is thus 
related to the continuous dynamical system associated with \eq{ev} in a
more sophisticated manner than usual:
The state space of the discrete system is
$(\R \times \Vs ) \times (\R \times \Vs )$ since only this set allows the 
leapfrog method to be defined as a map of states into states. Only a subset
of this state space (e. g. the one given by \eq{euler}) corresponds to
potential initial states of the continuous system \eq{ev}.

Notice the graphical manifestation of the leapfrog algorithm: The red dots 
make up the discrete trajectory. 
In all but the first and the last trajectory points
--- let an arbitrary such point be designated $(t_i,\psi_i)$ --- 
there is a short solid red line which indicates the slope required by the
differential equation at this position. 
To each such short solid line, there is a parallel broken line
which connects an already constructed trajectory point 
$(t_{i-1},\psi_{i-1})$ with the 
point $(t_{i+1},\psi_{i+1})$ to be constructed next.
\begin{figure}
\centering
\mbox{
\includegraphics[width=100mm]{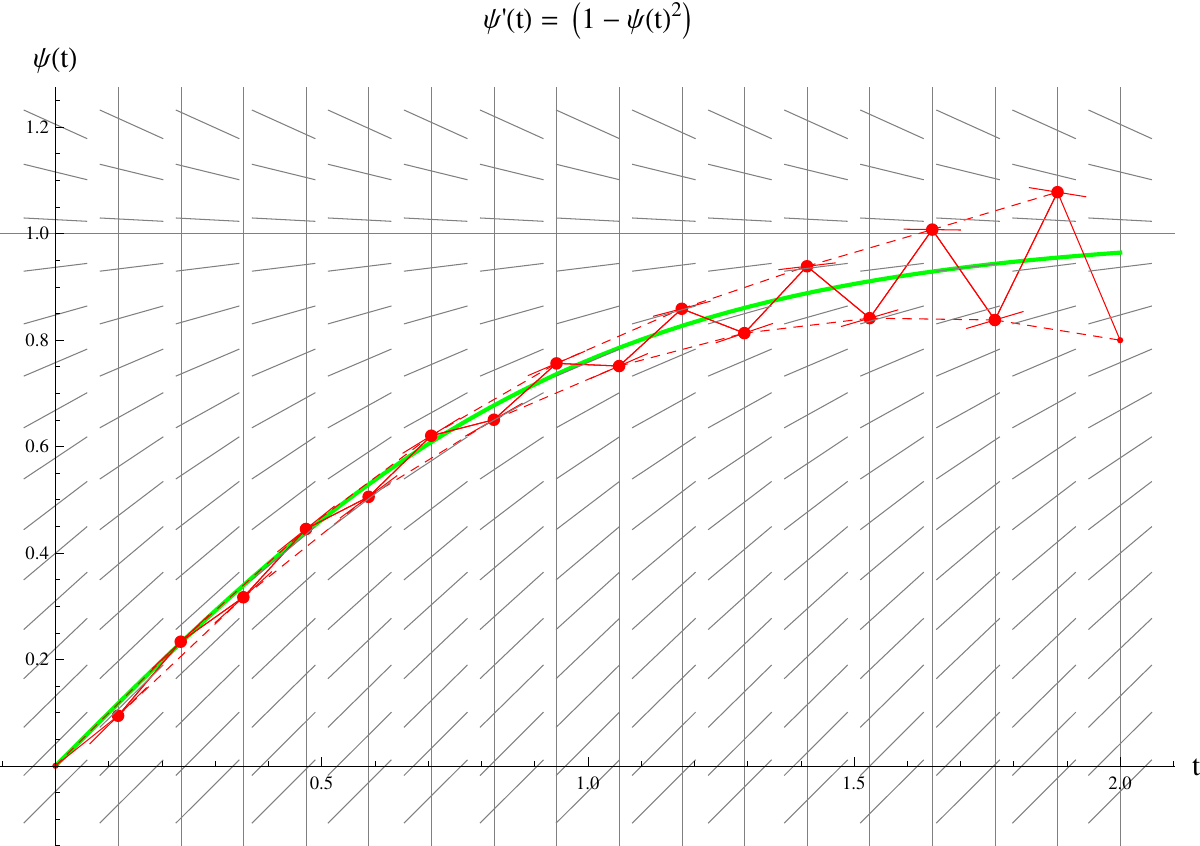}
}
\caption{Zigzag trajectory resulting from intentionally spoiled 
initialization.}
\label{fig:zigzag}
\end{figure}
The auxiliary lines, which are evident 
in Figure~\ref{fig:zigzag} are generated also for Figure~\ref{fig:Euler 
good}(b) but here they collapse here to a single polygon.

Equation \eq{iteration} defines $\psi_i$ for arbitrarily large $i$,
whereas equation 
\eq{ev} may drive a trajectory in finite time into infinity. 
A well-known example is 
\begin{equation}
\Lb{arctan}
\dot{\psi}(t) = 1 + \psi(t)^2 \ec \psi(0) = 0 \quad \text{thus} \quad 
\psi(t) = \tan(t) 
\quad \text{and} \quad \psi(\pi/2) = \infty \ep
\end{equation}
In such cases the size of the numbers involved in applying mapping $\LF$ 
repeatedly will grow above the size which can be handled with realistic
computational resources (which include computation time).
Such exploding situations also occur if the time step size is too large
for the differential equation under consideration.

It might be instructive to discuss the close correspondence of 
the leapfrog algorithm \eq{leapfrog} to the leapfrog game (Bockspringen in 
German). 
In the variant which is of interest here, there are two 
participants A and B in this nice dynamical sportive exercise. 
There is an intended direction of motion and the participants 
line up in this direction, B standing a few meters in front of A.
This is the initial condition, which corresponds in the algorithm to the 
ordered input $((t_0,\psi_0) \cong A, (t_1,\psi_1) \cong B)$. 
After two or three energetic steps, A jumps
over B, supporting himself with both palms on the shoulders of B, thereby 
receiving from 
B a smooth kick which makes A fly to a position sufficiently far in front 
of B that
now the action can continue with the roles of A and B reversed, then 
reversed again, and so forth. 
The kick which A receives from B corresponds to adding the term 
$(t_{2}-t_{0}) \,\dglx{\psi_1}{t_1}$
(associated with B) to the term $\psi_0$, which is associated with A.
The result of this addition is $\psi_2$, which corresponds again to A, but 
at a new position. 
By continuation we create terms $\psi_3, \psi_4 \ldots $~.
All terms with even index correspond to A, and those with odd index to B. 
The index grows with the progress along the intended direction of motion.

The algorithm can easily be shown to be \emph{reversible}: Let the
operator of motion reversal be defined as 
\fullmapeq{reversal}{\Ts}{(\R \times \Vs ) \times (\R \times \Vs )}{(\R 
\times \Vs ) 
\times (\R \times \Vs )}{((t_0,\psi_0), 
(t_1,\psi_1))}{((t_1,\psi_1),(t_0,\psi_0)) }
then we easily verify
\begin{equation}
\Lb{reversibility}
\LF \circ \Ts \circ \LF  = \Ts  \ec \Ts \circ \Ts = \id
\end{equation}
from which one concludes that $\LF$ is invertible, with the inverse given 
by
$\Ts \circ \LF \circ \Ts$.
This allows us to reconstruct from the last two data 
$((t_{n-1},\psi_{n-1}),(t_n,\psi_n))$
of a leapfrog trajectory all previous components 
$(t_k,\psi_k),\, k < n-1,$ 
\begin{equation}
\Lb{reversed}
(t_k,\psi_k) = \pi_{2}(\Ts \LF^{n-k} \Ts 
((t_{n-1},\psi_{n-1}),(t_n,\psi_n))) \ep
\end{equation}

If we would like to change time step size after having arrived at some
state $((t_{p-1},\psi_{p-1}),\,(t_{p},\psi_{p})$ to value $\tau$, we may
start a new synchronous trajectory with the state
\begin{equation}
\Lb{new step size}
((t_{p},\psi_{p}),\,(t_{p}+\tau,\psi_{p} + \tau F(t_p,\psi_p)) 
\end{equation}
or a potentially more accurate form which also involves $\psi_{p-1}$ 
(which equation \eq{new step size} simply forgets).

\section{An asynchronous version of the leapfrog method}
\label{asynchronous}

What I intend here, is to modify the leapfrog method in a way that no
trade-offs between simplicity and accuracy are involved when we start a
trajectory or when we change the time step size. 
A further aim is to 
preserve the computational simplicity of the algorithm.
In a narrower framework than \eq{ev} this modified leapfrog method has 
been introduced in \cite{qmi} and applied to time-dependent Hartree 
equations in \cite{sqd}.

We consider three consecutive components of a leapfrog trajectory of
\eq{ev}
\begin{equation}
\Lb{trajectory}
\stage{k} \ec \stage{k+1} \ec \stage{k+2} 
\end{equation}
and let $\tau$ be the time step: $\tau = t_{k+1} - t_k$.
Then we define velocity-like quantities $\phi$ as follows:
\begin{equation}
\Lb{phi}
\phi_k := \frac{\psi_{k+1} - \psi_k}{\tau} \ec \phi_{k+1} := 
\dgli{k+1} \ec \phi_{k+2} := \frac{\psi_{k+2} - \psi_{k+1}}{\tau} \ep
\end{equation}
From this definition and from \eq{leapfrog2} we obtain
\begin{equation}
\Lb{phi relation}
\frac{\phi_k + \phi_{k+2}}{2} = \frac{\psi_{k+2} - \psi_k}{2 \tau} = 
\dgli{k+1} = \phi_{k+1} \ep
\end{equation}
These equations allow us to compute $(t_{k+2},\,\psi_{k+2},\,\phi_{k+2})$ 
if
$(t_{k},\,\psi_{k},\,\phi_{k})$ and $\tau$ are given:
\begin{equation}
\Lb{algorithm}
\begin{split}
t_{k+1} &= t_k + \tau \ec \psi_{k+1} = \psi_k + \tau \,\phi_k \ec 
\phi_{k+1} = \dgli{k+1} \ec \\
t_{k+2} &= t_{k+1} + \tau \ec \\
\psi_{k+2} &= \psi_k + 2 \,\tau\, \phi_{k+1} \ec \\
\phi_{k+2} &= \frac{\psi_{k+2} - \psi_{k+1}}{\tau} \ep
\end{split}
\end{equation}
Equation \eq{phi relation} allows us to give the last two equations of 
\eq{algorithm} a more symmetrical form:
\begin{equation}
\Lb{algorithm second form}
\begin{split}
\phi_{k+2} &= 2\,\phi_{k+1} - \phi_k \ec \\
\psi_{k+2} &= \psi_{k+1} + \tau\, \phi_{k+2} \ep 
\end{split}
\end{equation}
The association 
$\mapto{(t_{k},\,\psi_{k},\,\phi_{k})}{(t_{k+2},\,\psi_{k+2},\,\phi_{k+2})}
$
can now be considered a mapping $\map{\R \times \Vs \times \Vs}{\R \times 
\Vs \times \Vs}$
which depends on the total time step $2\,\tau$ and thus will be denoted 
$\ALF_{2\tau}$.
The intermediary data $(t_{k+1},\,\psi_{k+1},\,\phi_{k+1})$ will be
referred to as the \emph{midpoint state} determined by
$(t_k,\,\psi_k,\,\phi_k)$ and $\ALF_{2\tau}$.
The mapping $\ALF_{2\tau}$ may be iterated to generate a discrete
trajectory $Tr$. 
This trajectory proceeds in time steps $2 \tau$, whereas the
leapfrog trajectory $Tr_0$ resulting from continuing \eq{trajectory}
proceeds in time steps $\tau$. $Tr$ thus sees $F$ only on a wider lattice
of $t$-values than $Tr_0$. So the two trajectories cannot be related in
any direct manner, although
$\ALF_{2\tau}$ was constructed out of two consecutive steps of $Tr_0$.
Nevertheless there is a leapfrog trajectory associated with $Tr$ in a
natural manner: it is given by the sequence of midpoint states of $Tr$.
Instead of iterating $\ALF_{2\tau}$ we could have considered a product
$\ALF_{2\tau_1}\circ\ALF_{2\tau_2}\circ\ALF_{2\tau_3} \ldots$ with varying
time step, and in this case no naturally associated leapfrog trajectory
can be expected to exist.

It is to be noted that going from the
normal leapfrog algorithm to the one given by the iteration of
$\ALF_{2\tau}$ we change the notion of system state.
The new state notion is more conventional in so far as it refers to a
single point in time, whereas the normal leapfrog state consists of data
that refer to two points in time.

\begin{figure}
\centering
\mbox{
\subfloat[Constant large time step.]
{\includegraphics[width=72mm]{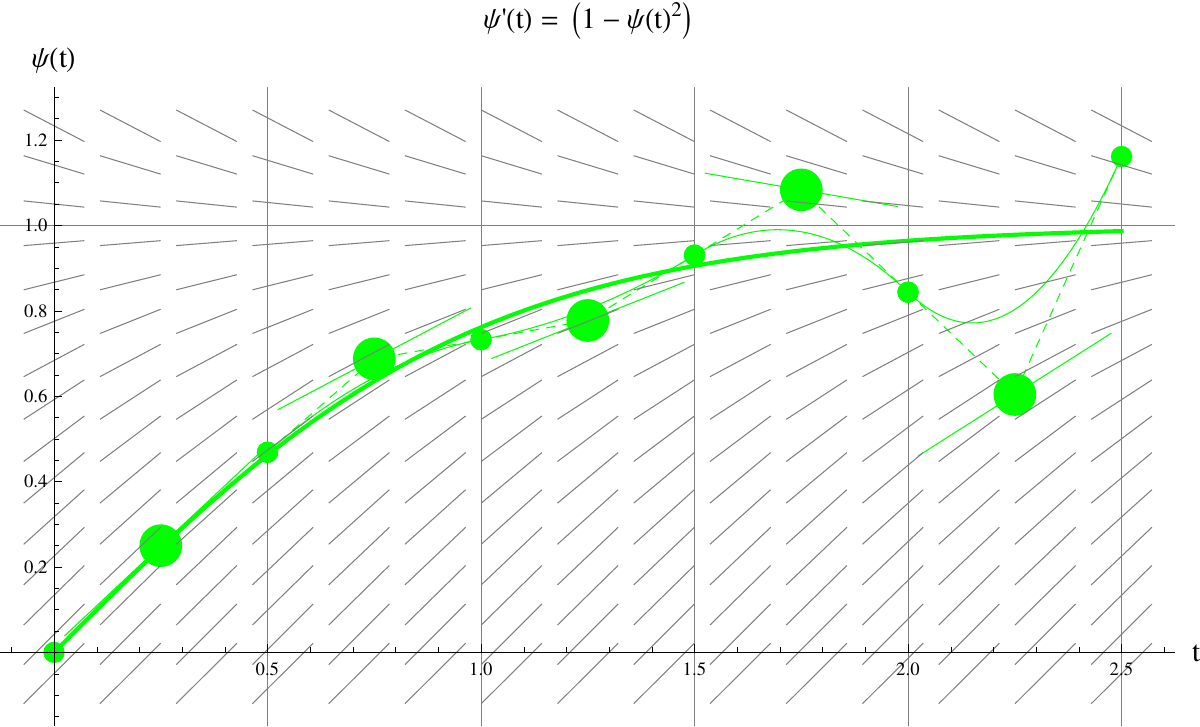}}
\quad
\subfloat[Half this time step in the second half of the trajectory.]
{\includegraphics[width=72mm]{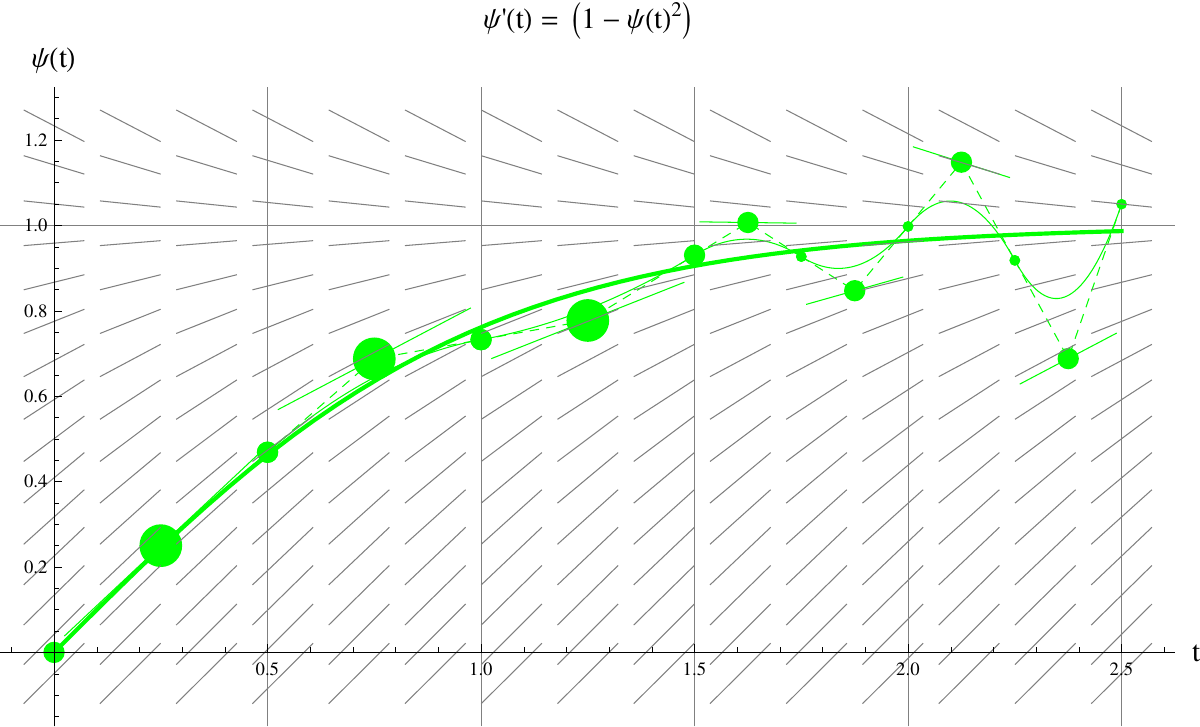}}
}
\caption{Geometry of \eq{algorithm third form} and reduction of the time 
step.}
\label{fig:ALF first}
\end{figure}

If we are given, according to the initial value problem of \eq{ev}, the 
initial values
$t_0$ and $\psi_0$, the augmentation to a full state according to the new 
state notion is straightforward and does not depend on the next time value
$t_1$. 
It is simply given by:
\begin{equation}
\Lb{initial phi}
\phi_0 := \dgli{0} \ep
\end{equation}
From $(t_0,\psi_0,\phi_0)$ and a not necessarily equi-spaced list
$(t_1,t_2,\ldots)$ we generate a discrete trajectory 
\begin{equation}
\Lb{trajectory0}
((t_0,\psi_0,\phi_0), (t_1,\psi_1,\phi_1), (t_2,\psi_2,\phi_2), \ldots) 
\end{equation}
by 
\begin{equation}
\Lb{iteration second}
(t_{i+1},\psi_{i+1},\phi_{i+1}) := \ALF_{h_i} (t_{i},\psi_{i},\phi_{i})
\ec 
h_i := t_{i+1}-t_i \ec
\end{equation}
where for each $h \in \R$ the mapping
\fullmapeq{Ah}{\ALF_h}{\R \times \Vs \times \Vs}{\R \times \Vs \times 
\Vs}{(t,\psi,\phi)}{(\u{t},\u{\psi},\u{\phi})}
is defined by the following chain of formulas:
\begin{equation} 
\Lb{algorithm third form}
\begin{split}
\tau &:= \frac{h}{2} \ec \\
t' &:= t + \tau \ec \\
\psi' &:= \psi + \tau \, \phi \ec \\
\phi' &:= \dglx{\psi'}{t'} \ec \\
\u{\phi} &:= 2\,\phi' - \phi \ec \\
\u{\psi} &:= \psi' + \tau\, \u{\phi} = \psi + 2 \tau \, \phi' \ec \\
\u{t} &:= t' + \tau \ep 
\end{split}
\end{equation}
Equations \eq{initial phi} to \eq{algorithm third form} define
the \emph{asynchronous leapfrog} method, sometimes abbreviated as
\emph{ALF}. 
This algorithm corresponds to 
equation (8) in \cite{qmi} but is more general since it does not assume the
special form of $F$ that was considered there.
Notice that the midpoint state data $t',\phi',\psi'$ appear only 
as intermediary quantities that help to give the algorithm an elegant form.
In particular, they do not belong to the discrete trajectory
\eq{trajectory0}, \eq{iteration second} generated by $\ALF$.
Their geometrical role becomes clear from the representation of the final 
state as
\begin{equation}
\Lb{inta}
\u{\psi} = \psi + h \,\phi + \frac{h^2}{2} \frac{\phi'-\phi}{\tau} \ec 
\u{\phi} = \phi + h \frac{\phi'-\phi}{\tau} \ep
\end{equation}
This representation suggests an interpretation in which $h$ is replaced by 
a parameter
which varies from $0$ to $h$ and thus connects the states $\psi$ and 
$\u{\psi}$ by a 
parabolic curve in the linear space $\Vs$ (and the quantities $\phi$ and 
$\u{\phi}$ by a linear curve). 
Everywhere along this connecting curve,
\emph{$\phi$ is the time derivative of $\psi$}.
The connecting parabola is easily seen to be the B\'{e}zier curve
generated by the
\emph{control points} $(t,\psi),(t',\psi'),(\u{t},\u{\psi})$.
In this way, the inherently time-discrete method proposes its own 
time-continuous representation. 
This is very convenient if one needs to compare 
trajectories from simulations with different time steps.
This time continuous representation is by mere interpolation; 
if one needs true detail about the history between $\psi$ and $\u{\psi}$ 
one has to reduce the time step in the simulation.
It is interesting to observe that the parabolas of adjacent time steps fit 
together
in a differentiable manner so that a sequence of time steps gives rise to
a quadratic  
\emph{B\'{e}zier spline} as a differentiable representation of the
discrete trajectory.
Figure~\ref{fig:ALF first} shows this spline curve together with the 
control points. 
The larger disks mark the intermediary configurations $(t',\psi')$ and the 
smaller
ones mark the configurations $(t,\psi)$ (or $(\u{t},\u{\psi})$) which 
belong to the discrete trajectory.
The short solid line attached to the larger disks indicates the direction 
given by the direction field of the differential equation. 
It coincides with the
direction determined by connecting the two neighboring smaller disks.
In Figure~\ref{fig:ALF first}(b), instead of the two final steps of 
sub-figure (a) we have four final steps of half the step size. Notice that
the size of the marking disks is coupled to the step size so that the
large disks belonging to the small steps equal in size just the small disks
belonging to the large steps. 
The disk at $t = 1.5$ marks the final discrete configuration
reached by a large step and also is the first discrete configuration from 
which a small step starts (so it has also to be marked with a disk half
this size;
the data structure of the graphics contains such a disk, it is hidden by 
the larger disk since it is not given a different color).

In the example of Figure~\ref{fig:ALF first} the horizontal course of the 
exact trajectory 
and the direction field provided by the differential equation in its 
neighborhood enforce the formation of a wave.
The significance of this phenomenon is not clear.
It is tempting to speculate that classical particle trajectories could be 
transformed to wave-like
processes by discretization. Some form of discretization should be
expected to happen, 
since the `computational resources of Nature' available for the evolution 
of any particle should be expected to be limited.

\begin{figure}
\centering
\mbox{
\includegraphics[width=140mm]{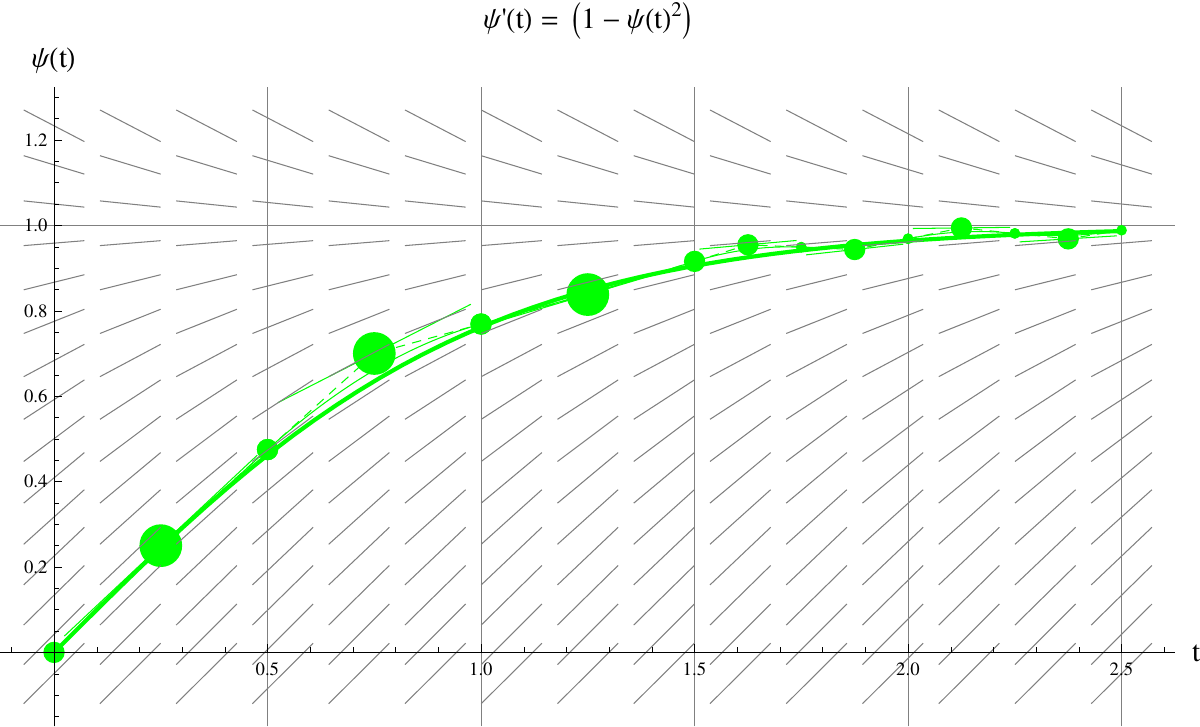}
}
\caption{Situation of Figure~\ref{fig:ALF first}(b) with relaxation 
($\lambda = 0.8$).}
\label{fig:relaxation}
\end{figure}

\begin{figure}
\centering
\mbox{
\subfloat[$\lambda = 0.95$]
{\includegraphics[width=70mm]{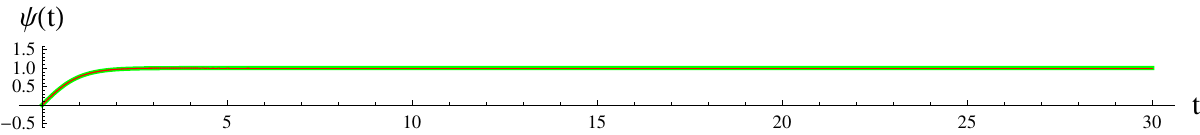}}
\quad
\subfloat[$\lambda = 0.975$]
{\includegraphics[width=70mm]{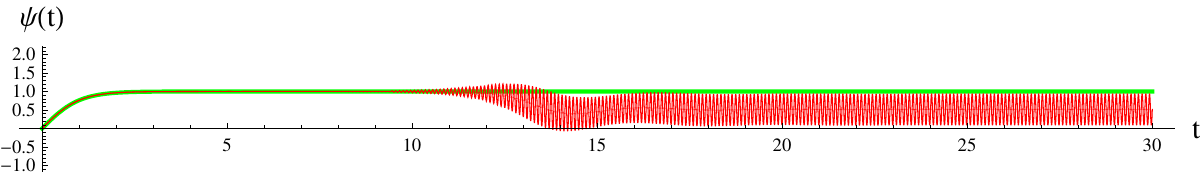}
}
}
\mbox{
\subfloat[$\lambda = 0.9875$]
{\includegraphics[width=70mm]{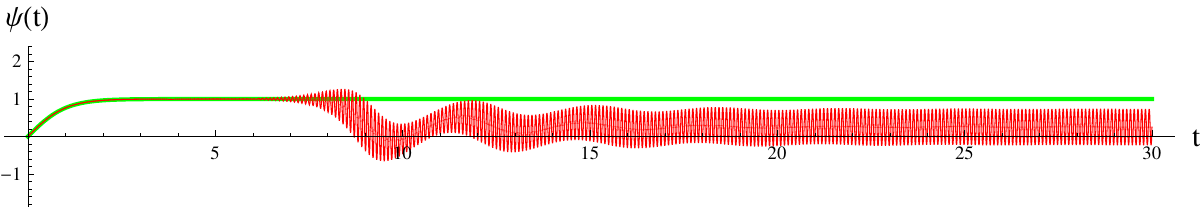}
}\quad
\subfloat[$\lambda = 1.0$]
{\includegraphics[width=70mm]{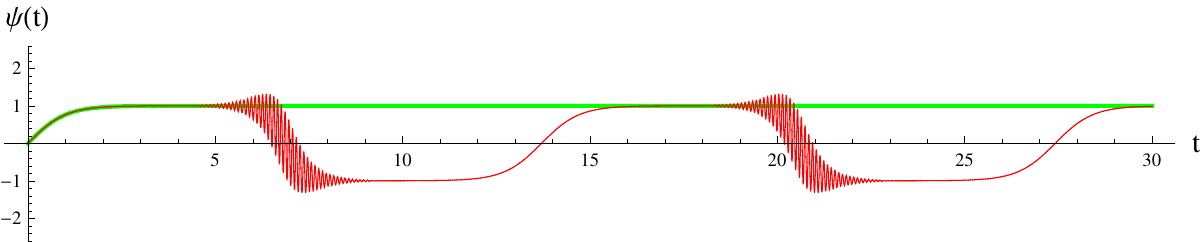}
}
}
\caption{How relaxation works for a virtually horizontal exact trajectory.}
\label{fig:ALF relaxation}
\end{figure}

The evolution equations \eq{algorithm third form} can be given a form
where no quantity needs to be copied and memorized:
\begin{equation} 
\Lb{algorithm third form revised}
\begin{split}
t \;&+=\; \tau \ec \\
\psi \;&+=\; \tau \, \phi \ec \\
\phi \;&+=\;2\, \lambda \,(\dglx{\psi}{t} - \phi) \ec \\
\psi \;&+=\; \tau\, \phi  \ec \\
t \;&+=\; \tau \ec 
\end{split}
\end{equation}
This property is obviously advantageous if $\psi$ and $\phi$ are large 
arrays of data
as they are in simulations of systems with many degrees of freedom. 
I tend to favor this property also at a conceptual physical level. 
The \emph{relaxation parameter} $\lambda$ introduced here has to be $1$ 
for \eq{algorithm third form revised} to be equivalent to \eq{algorithm 
third form}. 
Values slightly less than $1$ let the method work in some 
cases well where otherwise large deviations from the exact trajectory
would occur.
Figure~\ref{fig:relaxation} shows an example which demonstrates drastic 
reduction of the deviations seen in Figure~\ref{fig:ALF first}~. 
If we continue the trajectory to larger values of $t$, excessive 
oscillations seem to built up
unless relaxation is in place to prevent them. 
Figure~\ref{fig:ALF relaxation} illustrates this. 
Here the discrete trajectory is shown in red color by rendering only the
corresponding spline curve. 
The exact solution ($\tanh$, see \eq{tanh}) is shown in green color.
The phenomenon of oscillation is related to the property of reversibility: 
The sequence of $(\psi,\phi)$-states is not allowed to have two equal
components, since going back from two equal states a suitably
selected number of steps one would get two different points with $\psi=0$, 
although inspection of the discrete trajectory shows that there is only
\emph{one} such point. 
If the exact trajectory would not be effectively constant, the occurrence 
of equal $\psi$-values along the trajectory could easily be avoided. 
But in the case under consideration reversibility forces the
trajectory to make use of ever new values of $\psi$ and $\phi$ which is in 
conflict with the aim to render the exact trajectory with good accuracy. 
The mechanism which brings about this kind of `self-avoidance' seems to be 
poorly understood.

It might be convenient to see \eq{algorithm third form} rewritten in the 
self-explanatory style 
of \cite{Hairer}(which is similar to that in \cite{Reich}): 
We write our differential equation \eq{ev} as 
\begin{equation}
\dot{q} = f(t,q)
\end{equation}
and the initial data as $t_0, q_0$. We complement them by setting $v_0 := 
f(t_0, q_0)$
and are in a position to define a fully explicit time step which promotes 
data indexed by $n$ to
data indexed by $n+1$:
\begin{equation} 
\boxed{
\Lb{algorithm fourth form}
\begin{aligned}
t_{n+\frac{1}{2}} &= t_n + \frac{h}{2} \ec q_{n+\frac{1}{2}} = 
q_n + \frac{h}{2}\, v_n \ec v_{n+\frac{1}{2}} = f(t_{n+\frac{1}{2}}, 
\,q_{n+\frac{1}{2}})\ec\\
v_{n+1} &= 2\, v_{n+\frac{1}{2}} - v_n \ec q_{n+1} = 
q_{n+\frac{1}{2}} + \frac{h}{2}\, v_{n+1} = q_n + h\, v_{n+\frac{1}{2}}
\ec 
t_{n+1} = t_n + h \ep
\end{aligned}
}
\end{equation}
Notice that the initial condition plays a slightly exceptional role: For
$n 
> 0$ the quantities $v_n$ and $f(t_n,q_n)$
equal only approximately, whereas for $n = 0$ they equal exactly . 
The small quantity $\delta_n := v_n - f(t_n,q_n)$ is an interesting one to 
monitor in simulations. 
Notice $\delta_0 = 0$.

Now we return to our $\psi,\phi$-notation. 
The symmetric grouping of the formulas in \eq{algorithm third form} 
suggests that the time step map can be written
as a product of three maps.
For each $h \in \R$ we define the mapping
\fullmapeqnp{factor 1}{B_h}{\R \times \Vs \times \Vs}{\R \times \Vs \times 
\Vs}{(t,\psi,\phi)}{(t+h, \psi + h \, \phi, \phi)}
% \fullmapeq changed to \fullmapeqnp 2013-09-18
and the mapping
\fullmapeq{factor 2}{C}{\R \times \Vs \times \Vs}{\R \times \Vs \times 
\Vs}{(t,\psi,\phi)}{(t, \psi, 2\, F(t,\psi) -\phi)}
% \ep deleted 2013-09-18
We then see immediately
% than changed to then 2013-09-18
\begin{equation}
\Lb{factoring}
\ALF_h = B_{h/2} \circ C \circ B_{h/2} \ep
\end{equation}
Obviously $B_h$ is \emph{symplectic}: Writing $B_h(t,\psi,\phi)$ as 
$(t',\psi',\phi')$ 
(where now the apostrophe is used again as a normal diacritical mark and 
not in the special meaning of
\eq{algorithm third form})
we get
$\d \psi' \wedge \d \phi' = (\d \psi + h \, \d \phi) \wedge \d \phi = \d 
\psi \wedge \d \phi$.
Surprisingly $C$ is not symplectic but \emph{skew-symplectic}:
Writing $C(t,\psi,\phi)$ as $(t,\psi',\phi')$ we have 
$\d \psi' \wedge \d \phi' = \d \psi \wedge (2 \,\d F - \d \phi) = - \d
\psi 
\wedge \d \phi$,
since $\d F$ is proportional to $\d \psi$ (notice $\d t=0$) which implies 
$\d F \wedge \d \psi = 0$.
The product representation \eq{factoring} then implies that also $\ALF_h$ 
is 
skew-symplectic. Obviously the product of two skew-symplectic maps is 
symplectic.
\emph{Therefore, the densified form \eq{densified} of the asynchronous 
leapfrog integrator
is a symplectic explicit integrator.} Symplecticity plays here an unusual 
role, however,
since associating the `dynamical states' $(\psi,\phi)$ with a system 
trajectory is
a peculiarity of the method; the usual description uses $\psi$ alone.
In a conventional framework we encounter the variable $\phi$ only if our 
original
differential equation is the one which results from \eq{ev} by 
differentiation with 
respect to time:
\begin{equation}
\Lb{second order}
\begin{split}
\dot{\psi} &= \phi \ec \\ 
\dot{\phi} &= \frac{\partial F}{\partial t} + \frac{\partial F}{\partial 
\psi}\, \phi \ep
\end{split}
\end{equation}
Using integrator \eq{algorithm third form} for this equation may be 
non-trivial since
the integrator uses $F$ whereas the differential equation \eq{second
order} 
gives only $\frac{\partial F}{\partial t}$
and $\frac{\partial F}{\partial \psi}$, so that one obtains $F$ by solving 
a (probably partial)
differential equation.
The method to convert a first-order differential equation into a second 
order one by 
differentiation with respect to time and then to apply an explicit 
St\o rmer-Verlet integrator
is the path which led me to \eq{algorithm fourth form}. I followed this 
path in order 
to simulate quantum mechanical systems \cite{qmioriginal}, \cite{qmi}.

Obviously $\ALF_0$ is the identity map and $\ALF$ is \emph{reversible} in 
the sense 
that for all $h \in \R$ we have --- by a non-trivial cancellation of terms 
--- the equation
\begin{equation}
\Lb{reversible}
\ALF_{-h} \circ \ALF_h = \ALF_0
\end{equation}
which implies that each of the maps $\ALF_h$ is invertible, as is the 
product of 
arbitrarily many such mappings. As mentioned already for the corresponding
situation of the normal leapfrog method, this invertibility of all discrete
evolution maps does not imply that the dynamical system defined by 
\eq{ev} has invertible evolution maps.
The reversion of discrete trajectories is discussed in \cite{qmi} 
subsequent to
equation (10). It is to be noted that there is no motion reversion operator
comparable to \eq{reversal}. One may be tempted to try 
$\namedmapto{\Ts}{(t,\psi,\phi)}{(-t,\psi,-\phi)}$,
but this fails to satisfy $\ALF_h \circ \Ts \circ \ALF_h = \Ts$ which 
would 
correspond to \eq{reversibility}.
One should also note that the concept \eq{reversible} does not assume that 
the 
differential equation \eq{ev} satisfies any reversibility condition,
in particular not the one assumed in \cite{Holder}, after equation (1).

For a given state $(t,\psi,\phi)$, which determines a 
discrete trajectory by successive application of $\ALF_h$, one may
consider the leapfrog state $(t-h,\psi-h \phi),(t,\psi)$ which, by
successive application of $\LF$, creates a leapfrog trajectory the
first few steps of which agree reasonably well with those of the trajectory
considered before.
Therefore, in a sense, the role of $\phi$ is to memorize information from 
the foregoing integration step in addition to state data $\psi$. 
Also multi-step methods and predictor-corrector methods (or, more
generally, the \emph{general linear methods} of John Butcher \cite{glm},
which comprise all methods under current consideration)
improve computational economy by memorizing results from 
antecedent integration steps.
But they do so directly, by memorizing a list of previous states, each 
associated with the time of its validity.
Letting information from antecedent integration steps propagate in the form
of derived quantities, such as our $\phi$-data, was considered by me as a
unique characteristic of the present method. 
However, as J.M. Sanz-Serna pointed out to me, a similar velocity-related
quantity, often called the Nordsieck vector, is a constitutive element of
the Nordsieck method \cite{Nordsieck}, \cite{Osborne}~. 
As with our present method, it's role is to facilitate the change of the
step size. 
But unlike our $\phi$ the Nordsieck vector does not replace multi-step
values of the underlying multi-step method but serves as an additional
representation that is used in parallel with the multi-step values.
Nordsieck's method is intended to work for general multi-step methods and
therefore can't take advantage of the especially simple circumstances
holding for the leapfrog method.  

The reversibility of an integrator for the general equation \eq{ev}
implies more miracles than it implies for the reversible equations of 
Newtonian mechanics. 
This is pointed out in a discussion of the \emph{leaking bucket equation}
in \cite{bucket}~.
For this equation the final part of the exact trajectory is exactly 
horizontal so that we have similar but more transparent conditions as those
discussed above in connection with Figure~\ref{fig:ALF relaxation}~.

\section{The Kepler oscillator as a test example}
\label{Kepler}
Computing the motion of a point mass in the gravitation field of a 
stationary 
point mass is what the \emph{Kepler problem} is about. 
The radial motion in elliptic
Kepler orbits (as opposed to parabolic and hyperbolic ones) is oscillatory
and can be viewed as the motion of a one-dimensional oscillator which 
deserves interest as a mechanical example system. 
Unlike other non-linear example oscillators
such as the \emph{Duffing oscillator} and the \emph{Van der Pol oscillator}
this system seems to be anonymous. 
The self-suggesting name 
\emph{Kepler oscillator} can be found in \cite{Gallavotti} for this system
and will be used in the present article.
In the literature this name is, however, more often used for the harmonic 
oscillator 
which is related to the Kepler problem by a regularizing transformation,
known as the \emph{KS transformation}.
 
\begin{figure}
\centering
\mbox{
\includegraphics[width=125mm]{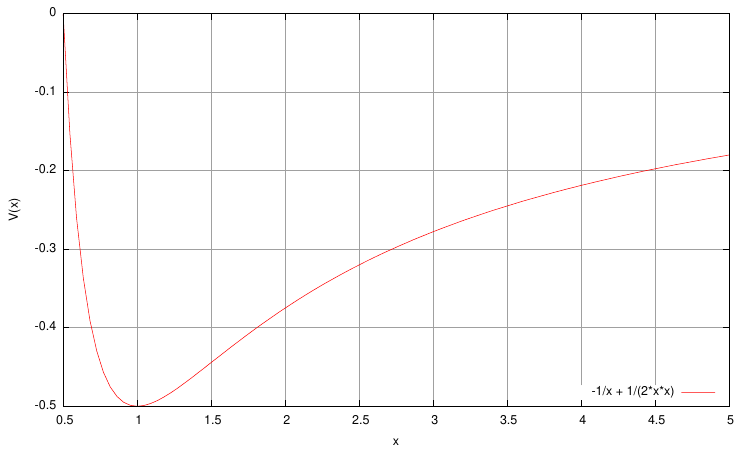}
}
\caption{Potential function of the Kepler oscillator.}
\label{fig:pubV}
\end{figure}

As is well known (e.g. \cite{Goldstein}, equation (3.14)) the radial
Kepler motion is governed by the differential equation
\begin{equation}
\Lb{radial}
m \ddot{r} = - \frac{\partial}{\partial r}\left( -\frac{G M m}{r} + 
\frac{L^2}{2 m r^2} \right)
\end{equation}
in which $L$ is the constant angular momentum of mass $m$ relative to the 
position of the
space-fixed mass $M$. Of course, $r$ is the distance between these two 
masses and $G$
is the constant of gravity. 
Restricting ourselves to orbits with non-vanishing $L$ 
and by selecting suitable units of time, mass, and length,
we get for the quantities $m$, $GM$, and $L$ the common numerical value 1. 
Writing $x$ for the numerical value of $r$ and $v$ for the numerical value 
of $\dot{r}$ we get 
\begin{equation}
\Lb{radial second}
\dot{x}=v \ec \dot{v} = - \frac{\partial}{\partial x}\left( -\frac{1}{x} + 
\frac{1}{2 x^2} \right) = \frac{1}{x^2}\left(\frac{1}{x}-1\right)
\end{equation}
which is the differential equation determined by \eq{radial} and also 
is (since, due to $m=1$, $v$ is the momentum) the system of canonical 
equations associated with the Hamiltonian
\begin{equation}
\Lb{Hamilton}
H(v,x) := T(v) + V(x) := \frac{1}{2} v^2 + \frac{1}{x} \left 
(\frac{1}{2x}-1 \right ) \ep
\end{equation}
This quantity is known to be constant on each orbit. Since, as 
Figure~\ref{fig:pubV} shows,
$V$ attains an absolute minimum at $x=1$: $V(1)=-\half $ we have $H(v,x) 
\geq H(0,1) = - \half$.
We consider only states for which $H(v,x)<0$ and thus $x>\frac{1}{2}$~. 
These correspond to the elliptical orbits in the Kepler problem; for them 
the radial motion has 
an oscillatory character. 
Kepler's ingenious method for computing the system path for given initial 
state
(not simply the orbit, a subject to which surprisingly many physics texts
restrict their interest) can be formulated as a simple algorithm: 
Given $(t_0,x_0,v_0)$ such that $H_0 :=H(v_0,x_0)<0$ and $t_1$ we have to 
go through the 
following chain of formulas (see also \cite{predicting}, Section 4):
\begin{equation}
\Lb{Kepler}
\begin{split}
   a&:= - \frac{1}{2 H_0} \text{  (major semi-axis) }\\
   \epsilon&:=\sqrt{1 - \frac{1}{a}} \text{  (numerical eccentricity) }\\
   n&:=a^{-\frac{3}{2}} \text{  (mean motion) }\\
   z&:= 1-\frac{x_0}{a} + \i \frac{x_0 v_0}{\sqrt{a}}\\
   E_0&:=\arg{z} \text{  (eccentric anomaly) }\\
   M_0&:=E_0 - \epsilon \sin{E_0} \text{  (mean anomaly) } \\
   M_1&:=M_0+(t_1-t_0)n \\
   E_1&:= \text{ solution of } E_1 = M_1 + \epsilon \sin E_1 \text{  
(Kepler's equation) }\\
   x_1&:=a(1-\epsilon \cos E_1) \\
   v_1&:=\frac{\epsilon a^2 n \sin E_1}{x_1}
\end{split}
\end{equation}
to get the exactly evolved state $(t_1,x_1,v_1)$. Here the solution $E$ of 
$E = M + \epsilon \sin E$ is given by the algorithm (C++ syntax, \verb+R+ 
is the 
type for representing real numbers, i.e. \verb+typedef double R;+)
\footnotesize
\begin{verbatim}
R solKepEqu(R M, R eps, R acc)
// M: mean anomaly, eps: numerical eccentricity, acc: accuracy e.g. 1e-8
{
   R xOld=M+1000, xNew=M;
   while ( abs(xOld-xNew) > acc ){
      xOld=xNew;
      R x1=M+eps*sin(xNew);
      R x2=M+eps*sin(x1);
      xNew=(x1+x2)*0.5;   // My standard provision against oscillations.
         // Works extremely well
   }
   return xNew;
}
\end{verbatim}
\normalsize
The computational burden for \eq{Kepler} is independent of the time span 
$t_1-t_0$,
and it does not matter whether this span is positive (prediction) or 
negative (retro-diction).
Hence, there is no relevant distinction between solution \eq{Kepler}
and what normally is referred to as a \emph{closed form solution}.
So, in assessing the accuracy of numerical integrators, we have the exact 
solution 
always available.
In addition to the original leapfrog method and the new asynchronous 
leapfrog
method, we consider two established second order methods for further 
comparison: 
The traditional 
\emph{second order Runge-Kutta method} (e.g. \cite{Press}, (16.1.2)) 
and the more modern symplectic \emph{position Verlet integrator}, 
\cite{Tuckerman}, equation (2.22). 
For this method there are several names in use, cf. \cite{qmi}, above 
equation (13), and \cite{Hairer}. 
The present article refers to it as the \emph{direct midpoint integrator} 
and recalls
its definition for the present simple situation that the forces don't 
depend on the velocity.
Equation \eq{radial second} can be viewed as a single differential
equation 
of second order
\begin{equation}
\Lb{differential equation of second order}
\ddot{x} = \frac{1}{x^2}\left(\frac{1}{x}-1\right)
\end{equation}
and for convenience of comparison with \eq{algorithm third form} we write 
this equation in a
form similar to \eq{ev} as
\begin{equation}
\Lb{ev second}
\ddot{\psi}(t) = \dglrhst \ec \dot{\psi}(t) =: \phi(t) \ep
\end{equation}
Since the differential equation is second order, the initial values for 
$\psi$ \emph{and} $\phi$
have to come from the problem and the integrator, just as in \eq{algorithm 
third form}, has the
task to promote them both. 
This is done by formulas very similar to \eq{algorithm third form}:
\begin{equation}  
\Lb{dmi}
\begin{split}
\tau &:= \frac{h}{2} \ec \\
t' &:= t + \tau \ec \\
\psi' &:= \psi + \tau \, \phi \ec \\
\u{\phi} &:= \phi + h\, F(t',\,\psi') \ec \\
\u{\psi} &:= \psi' + \tau\, \u{\phi} \ec \\
\u{t} &:= t' + \tau \ep 
\end{split}
\end{equation}
If one accepts to have one equation more than necessary, one may take 
\eq{algorithm third form}
as it stands, and replace the defining equation for $\phi'$ by the 
definition 
$\phi' := \phi + \tau \, \dglx{\psi'}{t'}$.

The implementation code for the integrators under consideration is 
contained in \verb+class KepOsc+
in file \verb+tut2.cpp+ which is listed in \cite{keposc}.

The orbits of the system can be uniquely parametrized by the values $0
\leq 
\epsilon < 1$ of the
numerical eccentricity, which is related to the total energy $H_0$ through
$\epsilon^2 = 1 + 2 H_0$~ (see \eq{Kepler}, notice that $H_0 < 0$). The
$x$-values along an orbit range
between the solutions $x_{\text{min}}, x_{\text{max}}$ of $V(x) = H_0$ and 
the $v$-values range
between the solutions $v_{\text{min}}, v_{\text{max}}$ of $T(v) - \half = 
H_0$.

\begin{figure}
\centering
\mbox{
\includegraphics[width=125mm]{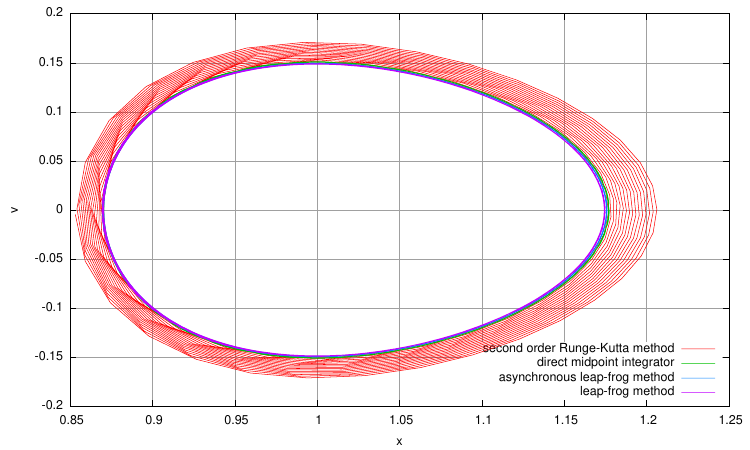}
}
\caption{Computed orbits in phase space.}
\label{fig:pub1a}
\end{figure}

\begin{figure}
\centering
\mbox{
\includegraphics[width=125mm]{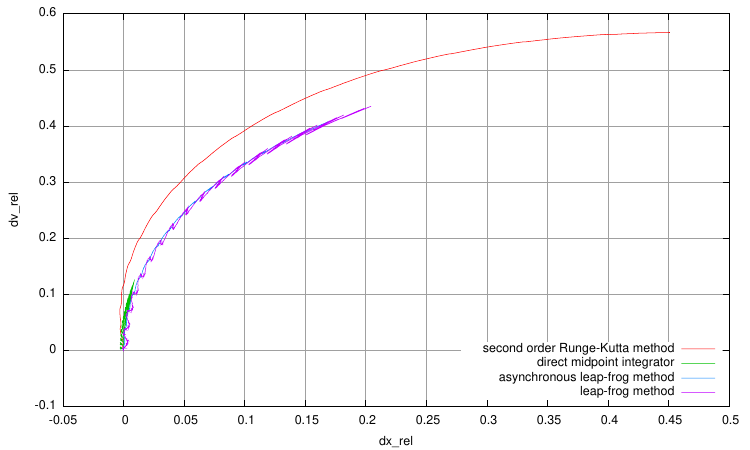}
}
\caption{Integration error of four computational methods in the 'numerical 
interaction picture'.}
\label{fig:pub2a}
\end{figure}

Most graphs to be presented here refer to a single path: The one which as 
an orbit is characterized
by $\epsilon = 0.15$, and the initial state of which is the `perihelion' 
i.e. $v_0=0$ and that $x_0=x_{\text{min}}$. 
The oscillation period of this path turns out to be $t_P=6.501$. 
Further we easily find $v_{\text{max}}=-v_{\text{min}}=0.157$
and $x_{\text{min}}=0.870$, $x_{\text{max}}=1.176$ which
agrees with the location of the oval shape in Figure~\ref{fig:pub1a}. 
As time step for stepwise integration we use $h=t_P/32$
to the effect that $32$ computed steps cover the whole period and thus 
would
lead back to the initial position if there would be no integration errors.
Each computation yields a discrete trajectory of $512$ steps, which 
corresponds to $16$ full
'revolutions'.
Figure~\ref{fig:pub1a} shows that for these data the Runge-Kutta method
does not create a periodic orbit and that the orbit in phase space
spirals into the outer space. 
At this graphical resolution, orbits created by the other methods are hard 
to distinguish.
Figure~\ref{fig:pub2a}
represents the deviation of the computed position from the exact one for
the four methods under consideration. 
What is displayed here is not simply the difference in phase space location
but the phase space position that occurs if the state at time $t$ is
back-evolved via the exact dynamics to the initial time $t=0$.
If the stepwise integration would not introduce an error, the point to be 
displayed would 
come out as $(0,0)$ in all cases.
The errors express themselves as curves (paths) with parameter $t$ and a 
longer curve
indicates a larger total error after the whole integration. 
Although the dependence on the curve parameter $t$ is not shown in the 
curves 
(only the orbit is represented),
the $16$ approximately repeated substructures in these curves show how the 
error evolves from revolution to revolution.
The coordinates in these diagrams are indexed `relative' which means that 
$x$-differences 
are divided by $x_{\text{max}} - x_{\text{min}}$
and $v$-differences are divided by $v_{\text{max}} - v_{\text{min}}$.
As already pointed out in \cite{svi}, near equation (92), these curves can 
be interpreted as paths in an \emph{interaction picture} dynamics, which
again is a dynamical system. 
This kind of interaction picture considers
the stepwise integration as the combined action of the exact evolution 
and a `discretization interaction' (in analogy to considering a digitized 
signal
as a superposition of the original analog signal and 'digitization noise').
It may therefore be fittingly referred to as \emph{numerical interaction 
picture}.
The more accurate the stepwise integration method, the weaker is the
interaction and the slower is the motion seen in the numerical interaction 
picture.
As mentioned in \cite{svi}, this diagnostics based on the numerical 
interaction picture 
is not restricted to systems for which the exact solution is directly 
accessible;
an access through stepwise back-evolution methods is sufficient if these 
are, say, two orders of 
magnitude more accurate than the method under investigation.

\begin{figure}
\centering
\mbox{
\includegraphics[width=125mm]{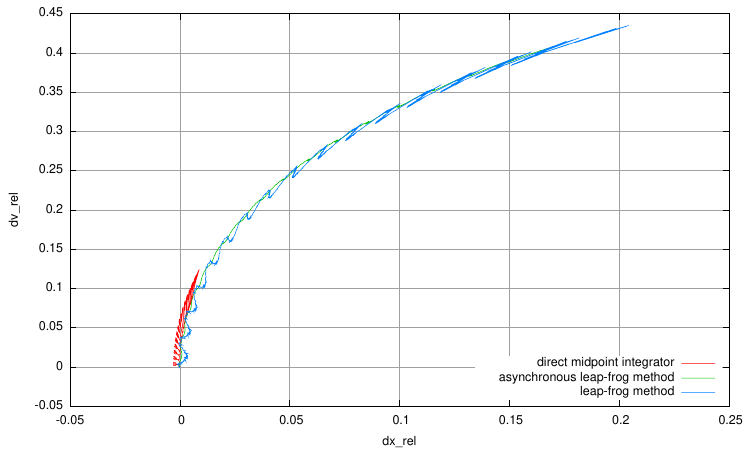}
}
\caption{Integration error of the better methods in the 'numerical 
interaction picture'.}
\label{fig:pub2b}
\end{figure}

\begin{figure}
\centering
\mbox{
\includegraphics[width=125mm]{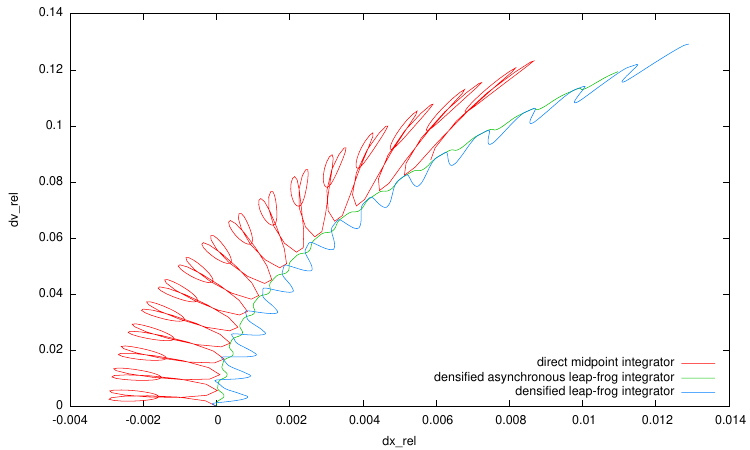}
}
\caption{Integration error of the direct midpoint method together with the 
densified
leapfrog methods in the 'numerical interaction picture' for 
$\epsilon=0.15$. Sixteen periods
at 32 steps per period.}
\label{fig:pub2c}
\end{figure}

\begin{figure}
\centering
\mbox{
\includegraphics[width=125mm]{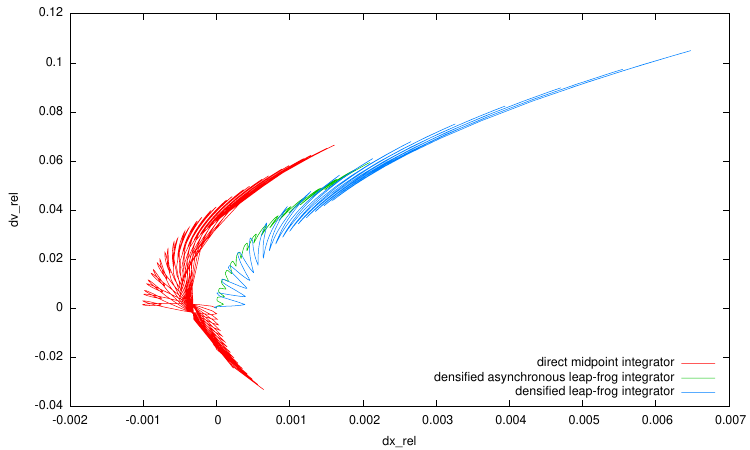}
}
\caption{Integration error of the direct midpoint method together with the 
densified
leapfrog methods in the 'numerical interaction picture' for 
$\epsilon=0.30$. Sixteen periods
at 64 steps per period.}
\label{fig:pub2d}
\end{figure}

This numerical interaction picture dynamics 
is related to the usage of evolution operators
$\e^{\i H t}\, \e^{-\i H_0 t}$ in quantum mechanical scattering theory and
with backward error analysis in numerical analysis of differential 
equations,
\cite{Calvo}, Chapter 10 , \cite{Hairer}, Section 4, \cite{Reich}, Chapter 
5.
The aim of \emph{backward error analysis} is to represent the 
discretization interaction
by additional terms to the right-hand side of the differential equation;
in case of Hamiltonian systems by an addition to the Hamiltonian, which is
the way for introducing interaction physicists are most familiar with.
The idea of the presently proposed method is to work with the dynamical 
systems 
directly without being forced to construct an equivalent Hamiltonian, or 
--- more generally --- a modified equation. The two movies \cite{movies} 
'Deformation of a phase space subset by interaction picture dynamics' show
features
\footnote{If one computes numerically the area enclosed by the curves
$\mapto{s}{(x(s),v(s))}$ 
which are shown in the movies, one finds them wiggling around a constant 
value for DALF-trajectories, and slowly linearly decreasing for 
ADALF-trajectories. Interestingly, for DALF-trajectories one finds
transient crossing over (and thus negative contributions to the enclosed
area) in these curves. Such a point of crossing over determines two
values $s_1$ and $s_2$ with $(x(s_1),v(s_1)) = (x(s_2),v(s_2))$,
but the full states $(x(s_1),v(s_1),w(s_1),a(s_1))$ and
$(x(s_2),v(s_2),w(s_2),a(s_2))$ (see \eq{adalf2 step} for the nomenclature)
do not agree to the effect that more than one DALF-trajectories may run
through one graphical point.}
which one would not easily read from any modified equation of 
backward error analysis (the converse is probably also true; it is the 
multitude of non-trivial observations which enhances our understanding).
Already our trajectory representation in 
Figures~\ref{fig:pub2a}, \ref{fig:pub2b}, \ref{fig:pub2c}, \ref{fig:pub2d} 
shows more morphologic features than the energy-error curves that
are normally generated as a kind of fingerprint of an integrator (e.g. 
\cite{Reich},
Fig 4.1).

Since, as already clear from Figure~\ref{fig:pub1a}, the error of the 
Runge-Kutta
method is much larger than the error of the other methods, the following 
Figure~\ref{fig:pub2b} gives the corresponding representation for 
the more accurate methods only.
This figure suggests that the direct midpoint integrator is by far more 
accurate
than the two leapfrog integrators. 
We will see now that this suggestion is misleading.
In the application considered in \cite{sqd} it was found that the 
asynchronous
leapfrog integrator showed similar step size requirements as the
direct midpoint integrator when a actual leapfrog step was defined as 
consisting of two leapfrog steps of half the step size. 
Also in the present context it makes sense to consider such a subdivision 
of a step.
We thus define the \emph{densified} leapfrog integrators as
\begin{equation}
\Lb{densified}
\begin{split}
\tilde{\LF} &:= \LF \circ \LF \\
\tilde{\ALF}_h &:= \ALF_{h/2} \circ \ALF_{h/2}
\end{split}
\end{equation}
and display the resulting error orbits in
Figure~\ref{fig:pub2c}. 
Note that again there are $32$ integration steps per orbit, but these are
made from substeps so that only every second computed step results in a 
graphical point. 
For such a combined step, the computational burden is
the same as for one second order Runge-Kutta step, but the accuracy is much
better than for Runge-Kutta.
It is plausible that only this densified version of the leapfrog methods 
turns
out to come close to the accuracy of the direct midpoint method: The latter
has direct access to the second derivative of the solution, whereas the 
leapfrog
methods only accesses the first derivative and thus can be viewed as 
simulating 
access to the second derivative by evaluating the first derivative at two 
different points.
One may generate corresponding diagrams for different values of 
eccentricity and
step size and will experience a surprising morphological stability of the 
curves and their relative length. 
Figure~\ref{fig:pub2d} is an example for this.
Here, the number of points per revolution is increased to 64 in response
to the increased value of the eccentricity $\epsilon$.

It is rather evident from all such graphs is that the asynchronous
leapfrog method has shorter and more regular error orbits than the
standard leapfrog method.
\clearpage

\section{The averaged densified asynchonous leapfrog integrator}
\label{adalf}
Let us introduce short names for the integrators under consideration: 
\begin{enumerate}
 \item \emph{ALF} asynchronous leapfrog, \eq{algorithm third 
form},\eq{algorithm fourth form}
 \item \emph{DALF} densified asynchronous leapfrog, \eq{densified}
% \eq{algorithm fourth form} replaced by \eq{densified} 2013-09-18 
 \item \emph{ADALF} averaged densified asynchronous leapfrog, to be 
defined in \eq{adalf step}
\end{enumerate}
The problem to be faced is evident from Figure \ref{fig:ALF first}. 
Here, consecutive ALF-steps end in directions that deviate 
considerably from the correct direction.
The deviation occurs in an alternating fashion. 
This suggests to define the final direction of a double step (thus of a 
single DALF-step) as the mean value out of the two ALF steps.
Let us first write the DALF-step in the most economic form, i.e. in the 
style of \eq{algorithm third form revised} with $\lambda = 1$: 
\begin{equation} 
\Lb{dalf step}
\begin{split}
t \;&+=\; \frac{\tau}{2} \ec \\
\psi \;&+=\; \frac{\tau}{2} \, \phi \ec \\
\phi \;&+=\;2\,(\dglx{\psi}{t} - \phi) \ec \\
\psi \;&+=\; \tau\, \phi  \ec \\
t \;&+=\; \tau \ec \\ 
\phi \;&+=\;2\,(\dglx{\psi}{t} - \phi) \ec \\
\psi \;&+=\; \frac{\tau}{2} \, \phi \ec \\
t \;&+=\; \frac{\tau}{2} \ep 
\end{split}
\end{equation}
Here the time step of the whole DALF-step is $2 \tau$. When considering an 
ALF-step as consisting of three sub steps, we have combined the last third 
of the first ALF-step with the first third of the second ALF-step and have 
saved a bit of work compared to a simple succession of two ALF-steps. 
This gives us the right to add a bit of computational complexity to 
obtain the definition of the ADALF-step: 
\begin{equation} 
\Lb{adalf step}
\begin{split}
t \;&+=\; \frac{\tau}{2} \ec \\
\psi \;&+=\; \frac{\tau}{2} \, \phi \ec \\
\phi \;&+=\;2\,(\dglx{\psi}{t} - \phi) \ec \\
\phi_1 &= \phi \ec \\
\psi \;&+=\; \tau\, \phi  \ec \\
t \;&+=\; \tau \ec \\ 
\phi \;&+=\;2\,(\dglx{\psi}{t} - \phi) \ec \\
\psi \;&+=\; \frac{\tau}{2} \, \phi \ec \\
\phi &= \frac{1}{2} (\phi + \phi_1) \ec \\
t \;&+=\; \frac{\tau}{2} \ec 
\end{split}
\end{equation}
as the most obvious realization of the averaging strategy indicated above. 
The general structure of the strategy is as follows: If an iteration is 
defined by application of a function $x_{i+1} := f(x_i)$, then a more 
robust and fast way to approach its iteration limit is by a modified 
iteration scheme $ x' := f(x_i)$, $x'' := f(x')$, $x_{i+1} := (x'+x'')/2$ 
which combines the original iteration scheme with a simple low-pass filter.
 
Of course, the algebraic nature of the 
quantities $x$ has to support formation of a mean value. 
This method occurred already in the code following \eq{Kepler} for solving 
Kepler's equation, and it worked perfectly in countless iteration 
applications for me.
As we will see, ADALF adds to the success story of this method. 
Very probably this method is known and honored in some scientific 
community and I would be very interested to get notice from anybody who 
knows it.

Finally it may be convenient to write the ADALF algorithm down for a 
second order system
\begin{equation}
 \ddot{x} = F(t,x,\dot{x}) 
\end{equation}
which applies 
directly to mechanical problems with velocity-dependent forces.
Transformed to first order it reads
\begin{equation}
 \dot{x} = v \ec  \dot{v} = F(t,x,v) 
\end{equation}
with initial values $x(t_0) = x_0$, $v(t_0) = v_0$. 
We write $(x,v)$ for $\psi$ and $(w,a)$ for $\phi$ and translate 
\eq{adalf step} to
\begin{equation} 
\Lb{adalf2 step}
\begin{split}
t \;&+=\; \frac{\tau}{2} \ec \\
x \;&+=\; \frac{\tau}{2} \, w \ec v \;+=\; \frac{\tau}{2} \, a \ec \\
w \;&+=\;2\,(v - w) \ec a \;+=\;2\,(F(t,x,v) - a) \ec \\
w_1 &= w \ec a_1 = a \ec \\
x \;&+=\; \tau \, w \ec v \;+=\; \tau \, a \ec \\
t \;&+=\; \tau \ec \\ 
w \;&+=\;2\,(v - w) \ec a \;+=\;2\,(F(t,x,v) - a) \ec \\
x \;&+=\; \frac{\tau}{2} \, w \ec v \;+=\; \frac{\tau}{2} \, a \ec \\
w &= \frac{1}{2} (w + w_1) \ec a = \frac{1}{2} (a + a_1) \ec\ \\
t \;&+=\; \frac{\tau}{2} \ep 
\end{split}
\end{equation}
The initialization of $(x,v),(w,a)$ is of course given by $x(t_0) = x_0$, 
$v(t_0) = v_0$, $w(t_0) = v_0$, and $a(t_0) = F(t_0,x_0,v_0)$. The versions
without an averaging step (i.e. DALF) and without a separation in two 
half-steps (i.e. ALF) are obvious from this pattern.
Notice that in \eq{ev second}, \eq{dmi} we have formulated the the 
St\o rmer-Verlet integrator only for forces which don't depend on the
velocity, whereas here such forces are nothing particular.

\subsection{Stability}
\label{stability}
As is well known the behavior of integrator-created trajectories of the 
simple linear test equation 
\begin{equation}
\Lb{linear test}
\dot{\psi}(t) = \omega \,\psi(t)\ec \omega \in \C \ec |\omega| = 1
\end{equation}
for a complex-valued function $\psi$ helps to analyze the complex 
phenomenon of \emph{stability} of integrators \cite{Stuart}, Def 3.6.1, 
\cite{shampine}.
Due to the simple nature of \eq{linear test} the integrators of present 
interest can be represented by complex propagation matrices: 
\begin{equation}
 \Lb{matrix 1}
 \begin{pmatrix} \psi_{n+1} \\ \phi_{n+1} \end{pmatrix} = 
 \begin{pmatrix} \alpha(h, \omega) & \beta(h, \omega) \\ 
\gamma(h, \omega) 
& \delta(h, \omega)  \end{pmatrix} \begin{pmatrix} \psi_n \\ \phi_n 
\end{pmatrix} \ec
\end{equation}
where $h$ is the time step.
It is a straightforward, though tedious and error-prone task (thus best 
done with help from a computer algebra system), to compute these 
propagation matrices from the defining equations for the integrators. Here 
are the results:
The matrix elements $\alpha, \beta, \gamma, \delta$ are given for ALF:
\begin{alignat}{2}
 \Lb{matrix alf}
 \alpha(h, \omega) &= 1 + h \omega & \qquad 
 \beta(h, \omega) &= \frac{h^2 \omega}{2} \\
 \gamma(h, \omega) &= 2 \omega & \qquad 
 \delta(h, \omega) &= -1 + h \omega \ec
\end{alignat}
and for ALF with relaxation parameter $\lambda$
\begin{alignat}{2}
 \Lb{matrix ralf}
 \alpha(h, \omega) &= 1 + h \lambda \omega & \qquad 
 \beta(h, \omega) &= h \Bigl( 1 + \lambda \bigl( \frac{h \omega}{2} - 1 
\bigr) \Bigr) \\
 \gamma(h, \omega) &= 2 \lambda \omega & \qquad 
 \delta(h, \omega) &=  1 + 2 \lambda \bigl( \frac{h \omega}{2} - 1 \bigr) 
\ec
\end{alignat}
and for DALF:
\begin{alignat}{2}
 \Lb{matrix dalf}
 \alpha(h, \omega) & = 1 + h \omega + \frac{h^2 \omega^2}{2} & \qquad 
 \beta(h, \omega) & = \frac{h^3 \omega^2}{8} \\ 
 \gamma(h, \omega) & = 2 h \omega^2  & \qquad
 \delta(h, \omega) & = 1 - h \omega +\frac{h^2 \omega^2}{2} \ec
\end{alignat}
and for ADALF:
\begin{alignat}{2}
 \Lb{matrix adalf}
 \alpha(h, \omega) & = 1 + h \omega + \frac{h^2 \omega^2}{2} & \qquad 
 \beta(h, \omega) & = \frac{h^2 \omega (1 + h \omega)}{16} \\ 
 \gamma(h, \omega) & = 2 h \omega^2  & \qquad
 \delta(h, \omega) & = \frac{h \omega (-1 + h \omega)}{4} \ep
\end{alignat}
As will be seen in Subsection~\ref{order} these leapfrog methods 
are all of second order.
Therefore it is natural to compare them with the general second order 
explicit Runge-Kutta method. 
This method depends on a parameter $a_1$. 
Giving this parameter the value $\,0, 1/3, 1/2\,$ defines the midpoint 
method, Ralston's method, and Heun's method respectively. 
The only Runge-Kutta method considered in Section~\ref{Kepler} is the 
midpoint method.

In order to allow Runge-Kutta methods to be used 
interchangeably with leapfrog methods, especially in an analysis of 
their stability, we have to augment their state space by a velocity. 
This velocity is initialized at the beginning of 
a trajectory by the same formula as in the leapfrog case, and the 
evolution step is reformulated in a way that it takes the velocity as 
given at the beginning (saving one evaluation of the the right-hand side 
of the differential equation) and has to be 
updated at the end of the step (consuming one evaluation of that 
right-hand side). 
It is clear that this does not influence the states along a trajectory in 
any way, with the exception of the last state: In a trajectory of the 
modified style it determines not only a position but also a velocity. 
The computational work that generated this additional information was done 
with the initialization step. 
As we will see, the circumstance that we 
have a value for the velocity at the beginning and at the end of an
evolution step allows us to implement a cheap and effective automatic step 
control.
With the modification described here, this auto-step method works also 
with Runge-Kutta and an instructive direct comparison becomes possible.

The effect of this modification that is of main interest in the present 
context is that it allows us to describe the Runge-Kutta integration 
step by a propagation matrix just as in the previous cases: 
\begin{alignat}{2}
 \Lb{matrix runge kutta}
 \alpha(h, \omega, a_1) & = 1 + h \omega (1 - a_1)  & \qquad 
 \beta(h, \omega, a_1) & = h \bigl(a_1 + \frac{h \omega}{2}\bigr) \\ 
 \gamma(h, \omega, a_1) & = \omega (1 + h \omega (1 - a_1))  & \qquad
 \delta(h, \omega, a_1) & = h \omega \bigl(a_1 + \frac{h \omega}{2}\bigr) 
\ep
\end{alignat}

Successive application of integration steps translates to forming powers
of the propagation matrices. 
Stability of the integration method is related to the question whether 
powers remain bounded when the number of factors tends to infinity
\footnote{
Here we argue in the framework defined by \eq{linear test}. 
If $\omega$ has a positive real part, the exact solution grows 
exponentially. So,one would not expect $h \omega$ to belong to $\mathbf{S}$ 
for any $h > 0$, or --- put differently --- one expects that $\mathbf{S}$ 
has no part which belongs to $\set{z \in \C}{\Re(z) > 0}$. 
If in such a case the powers of the propagation matrix remain bounded this 
is a deficiency of the integration method.
}~. 
Fortunately all complex 2 by 2 matrices are diagonalizable and thus allow 
a simple representation of their powers. We thus easily see:
Let $A$ be any complex 2 by 2 matrix. Then $\text{sup}\: 
\set{\norm{A^n}}{n \in \N} < \infty$ iff (i.e. if and only if) 
$|\lambda_1| \le 1$ and 
$|\lambda_2| \le 1$, where $\lambda_1, \lambda_2$ are the eigenvalues of 
$A$.

The complex eigenvalues of the propagation matrices of 
the various methods can be computed straightforwardly: We have for ALF:
\begin{equation}
\Lb{ev alf}
 \lambda_{1,2}(h \omega) = h \omega \pm \sqrt{1 + h^2 \omega^2} \ec
\end{equation}
and for ALF with relaxation parameter $\lambda$
\begin{equation}
\Lb{ev ralf}
 \lambda_{1,2}(h \omega) = 
 1 + \lambda ( h \omega - 1 ) \pm 
 \sqrt{\lambda \bigl(2 h \omega + \lambda (h \omega - 1)^2 \bigr)} \ec
\end{equation}
and for DALF:
\begin{equation}
\Lb{ev dalf}
 \lambda_{1,2}(h \omega) = \frac{1}{2} \left( 2 + h^2 \omega^2 \pm    
h \omega \sqrt{4 + h^2 \omega^2} \right)\ec
\end{equation}
and for ADALF:
\begin{equation}
\Lb{ev adalf}
 \lambda_{1,2}(h \omega) = \frac{1}{8} \left( 4 + 3 h \omega + 3 h^2 
\omega^2 \pm    
\sqrt{16 h \omega + \left(4 + 3 h \omega + 3 h^2 \omega^2 \right)^2} 
\right)
\end{equation}
and for second order Runge-Kutta:
\begin{equation}
\Lb{ev rk}
\lambda_1(h \omega) = 0 \ec \lambda_2(h \omega) = 1 + h \omega + \frac{h^2 
\omega^2}{2} \ep
\end{equation}
Notice that the eigenvalues for Runge-Kutta are independent of the 
parameter $a_1$.

The number $h \omega \in \C$ belongs to the \emph{set of absolute 
stability} $\mathbf{S}$ (\emph{stability region} for short) of an 
integrator of the kind considered here, iff 
$|\lambda_1( h \omega)| \le 1$ and $|\lambda_2( h \omega)| \le 1$.
A value of $h$ such that $\omega h$ lies on the boundary of $\mathbf{S}$
is said to be \emph{critical} since for time steps larger than this $h$ 
the integration method will yield exploding trajectories. 
Time steps lower than a critical one will be called \emph{sub-critical}. 
As we have seen above, a time step is sub-critical iff there exists 
a finite upper bound for the set of all powers of the propagation 
matrix.
The qualification `absolute' refers to the fact that it is valid for all 
initial states and thus for all trajectories.
To see this, we consider the case that only one of the eigenvalues has an 
absolute value larger than one. 
The best chance to get a non-growing trajectory would be given if the
initial condition $(\psi_0,\omega \psi_0)$, in an expansion in terms of 
the two eigenvectors would have only one component, namely the one 
belonging to the lower eigenvalue. (This may not be possible since
the possible initial conditions form a linear subspace of our 
$(\psi,\phi)$-space.) Even in this best conceivable case a numerically
generated trajectory will explode in the long run since numerical noise 
will create admixtures of the second eigenvector and these will grow 
exponentially.

One should be aware this conventional notion of a `set of absolute 
stability' is a bit contrived in that it mixes quantities referring to 
different differential equations (namely different values of $\omega$) 
with quantities referring to different discrete trajectories of a single 
equation (namely different values of $h$). 
Here I use this notion only in order to enable comparison with the 
literature.  
The notion of a critical step size makes sense also for a 
single differential equation. 

\begin{figure}
\centering
\mbox{
\subfloat[Total view of the stability regions.]
{\includegraphics[width=72mm]{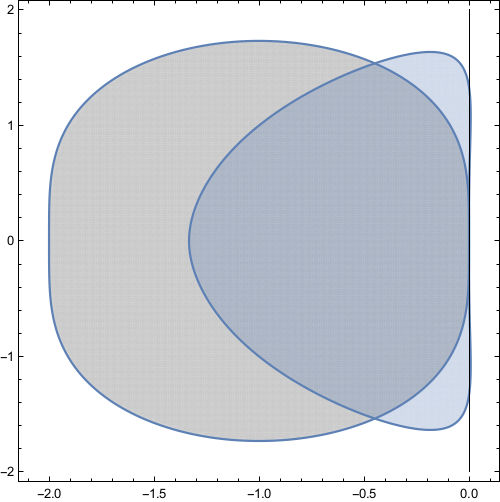}}
\quad
\subfloat[Detail view of the stability regions.]
{\includegraphics[width=72mm]{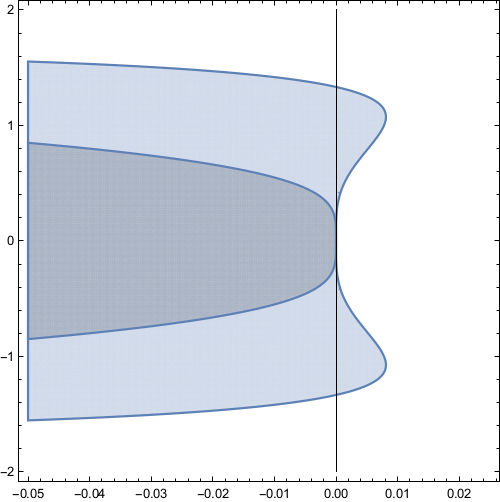}}
}
\caption{Regions of stability for 2nd order Runge-Kutta (gray), ADALF
(light blue), and DALF (a black line).}
\label{fig:pubxx}
\end{figure}

\begin{figure}
\centering
\mbox{
\subfloat[Total view of the stability regions.]
{\includegraphics[width=72mm]{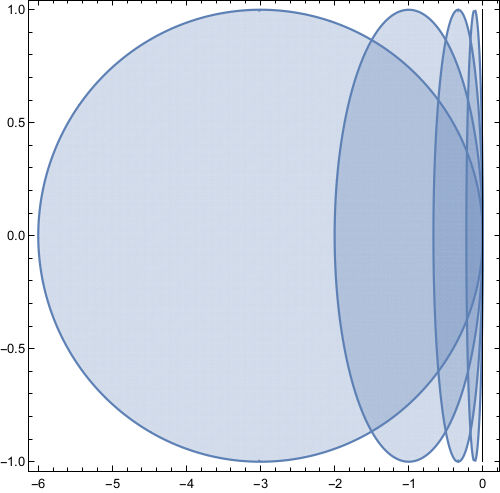}}
\quad
\subfloat[Detail view of the stability regions.]
{\includegraphics[width=72mm]{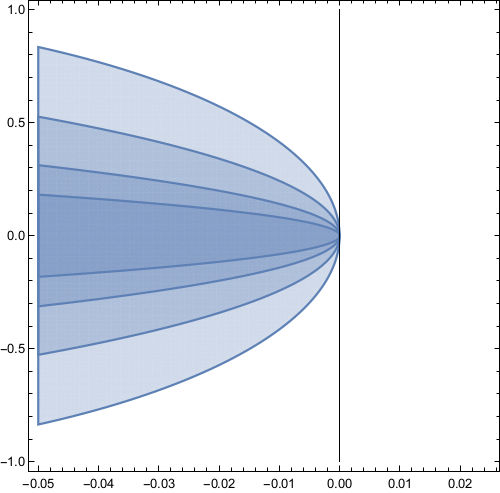}}
}
\caption{Regions of stability for ALF with relaxation $\lambda = 0.25, 
0.5, 0.75,0.9$ (light blue, from left to right), and without relaxation 
(a black line).}
\label{fig:stability-relaxed}
\end{figure}

Given the explicit formulas for $\lambda_1$, $\lambda_2$ it is 
straightforward to work out the regions of stability as done here by a 
commercial graphics function
\footnote{RegionPlot of Mathematica}
that created Figures~\ref{fig:pubxx} and \ref{fig:stability-relaxed}.

An instructive special case of our test equation \eq{linear test} is the 
Schr\"{o}dinger equation in one-dimensional Hilbert space. 
This is the case that $\omega$ is purely imaginary and the exact trajectory
is oscillatory (actually, $t \mapsto \exp(it)$). 
Suggesting a quantum mechanical context we write 
\[ \omega =  - \i H \quad \text{and} \quad \i \dot{\psi} = H \psi \ec\]
where $H \in \R$. 
The critical time step results from the section of 
$\mathbf{S}$ with the imaginary axis:
\begin{equation}
\Lb{stability size}
h_{\text{crit}}(\mathbf{S}) := \text{sup}\: \set{ r > 0}{ r \i \in 
\mathbf{S}} \ep
\end{equation}
As Figures~\ref{fig:pubxx} and \ref{fig:stability-relaxed} suggest, and is 
easily derived from the defining equations \eq{ev alf}, 
\eq{ev ralf}, \eq{ev dalf}, \eq{ev adalf}, \eq{ev rk},
this quantity has the value $0$ for the 2nd order Runge-Kutta methods and
for ALF with relaxation $\lambda < 1$. It has the value $1$ for ALF,
the value $4/3$ for ADALF, and $2$ for DALF. 

Let us put these numbers into perspective. We begin with method ALF.
The critical step size is $1$ which implies that there are $2 \pi$ steps 
per oscillation period. 
For method DALF the critical step size twice as large.
Since each step does 2 evaluations of the right-hand side of the 
differential equation, we get $2 \pi$ such evaluations per oscillation 
period --- just as much as with ALF.
Since the critical time step for ADALF is by a factor 1.5 ($= 2 : 4/3$)
shorter than the DALF step, we need $3 \pi$ evaluations per period. 
This roughly places the evaluations at the corners of a regular hexagon or 
a regular nonagon. 
It is plausible that for less evaluations it is hard to follow a circle.

That there is no sub-critical time step for the Runge-Kutta methods implies
that $\psi$ grows over any boundary in the long run. This is in remarkable 
contrast to the nice behavior of the leapfrog solutions for which --- in 
quantum mechanical terminology --- both the norm of the state and the 
expectation value of the energy are constant up to an additive fluctuation 
of order $O(h^2)$.
In the case of the DALF method this can be inspected in detail with the 
interactive application \cite{stability}.
This determination of critical step size generalizes naturally to quantum 
mechanics in finite-dimensional Hilbert spaces where $H$ is a large 
Hermitian matrix, see e.g. \cite{qmioriginal}, \cite{qmi}. For $\Delta 
t_{\text{crit}}\norm{H}$ one gets the same values which were given 
above for $h_{\text{crit}}$.
There are fast and cheap methods (e.g. 
\cite{qmioriginal} equ. 80) to compute the operator norm $\norm{H}$ for 
large matrices $H$.
An interactive DALF-simulation of a wave function evolving according to 
the Dirac equation in one-dimensional space under the influence of an 
adjustable electric potential is available at  \cite{klein paradox}. 
Here one can change the times step within the range of sub-critical values 
while the program is running and is showing the moving wave function on 
screen.
As a rule of thump it can be stated that time steps five to ten times 
smaller than the critical value let dynamical features appear virtually 
independent of the time step.
If one changes the time step to negative values one goes back to the 
initial condition with impressive accuracy.

There is a striking difference between the pure leapfrog methods, which 
gives $\mathbf{S}$ only an interval on the imaginary axis, and the 
averaged version, which gives a region that extends into the 
left half-plane and also, to a much smaller extend, to the right half-plane 
(see Figure~\ref{fig:pubxx}(b)). 

To explain this difference we consider the solutions of \eq{linear 
test} with initial condition $\psi(0) = 1$ for various values of $\omega$ 
(recall, $\omega \in \C, |\omega| = 1$). 
For $\omega = \i$ we have the purely oscillatory solution 
$\psi(t) = \e^{\i t}$ for which leapfrog works well for sub-critical 
values of the time step.
If, however, $\omega$ gets a non-vanishing negative real part, the exact 
solution is a damped oscillation that finally dies out and becomes 
horizontal and nearly straight. 
This, together with the reversibility of the leapfrog method, causes 
problems that motivated the relaxation parameter in 
\eq{algorithm third form revised} and were discussed together with the 
relaxation approach. 
This approach is now considered obsolete and the ADALF method is the 
proper replacement. 
As Figure~\ref{fig:pubxx}(a) indicates, both oscillatory and damped 
trajectories can be represented with ADALF using rather large time steps.
The common feature of the relaxation method and the ADALF method is that 
they break reversibility. 

\section{The Kepler oscillator as a test example, continued}
\label{Kepler 2}
\begin{figure}
\centering
\mbox{
\includegraphics[width=140mm]{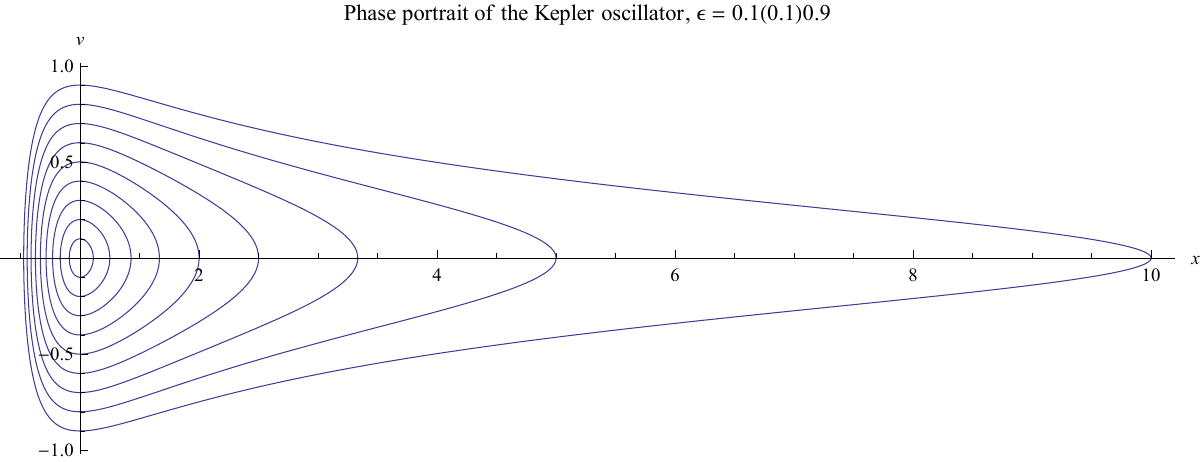}
}
\caption{Phase space portrait of the Kepler oscillator}
\label{fig:portrait}
\end{figure}

Compared to the treatment in Section~\ref{Kepler} we consider in some 
cases also trajectories with numerical eccentricity up to 
$0.99$. Figure~\ref{fig:portrait} shows the phase space portrait up to 
eccentricities $\epsilon = 0.9$. 
As is obvious from this portrait 
and even more so from considering the corresponding Kepler orbit,
for strongly eccentric orbits the velocity varies over a wide range so 
that we then have an example of a \emph{stiff system}. 
When we give for two time-stepped trajectories the same value for steps 
per revolution, although the eccentricities differ considerably, this is 
not a fair comparison.
Since the velocities vary more strongly for the more eccentric trajectory, 
one should use more steps for the more eccentric trajectory, if 
generating a discrete approximation with fixed time step should be an 
equally demanding task for both trajectories. 
To find the proper enhancement factor it is helpful to consider the 
full Kepler trajectory from which one gets our oscillator trajectory as
the projection onto the radial direction. 
Using the distance $r$ from the central body we have for the angular 
velocity $n$ (see \eq{Kepler}) the formula $L = r^2 n = 1$ (known as 
Kepler's second law, see \eq{radial} and following for $L$). From 
$r_\mathrm{min} = x_\mathrm{min} = 1/(1 + \epsilon)$ we get 
$n_\mathrm{max} = r_\mathrm{min}^{-2} = (1 + \epsilon)^2$ and 
\[n_\mathrm{mean} = a^{-3/2} = (1 - \epsilon^2)^{3/2} \ep\]
Then, the appropriate number of steps per revolution depends on the 
eccentricity as follows: 
\[ N_\mathrm{perRev}(\epsilon_1) : N_\mathrm{perRev}(\epsilon_2) = 
 \frac{n_\mathrm{max}(\epsilon_1)}{n_\mathrm{mean}(\epsilon_1)} : 
\frac{n_\mathrm{max}(\epsilon_2)}{n_\mathrm{mean}(\epsilon_2)} \ep\]
For $\epsilon_1 := 0$ and $\epsilon_2 := \epsilon$ we get 
\begin{equation}
\Lb{epsnorm}
 N_\mathrm{perRev}(\epsilon) = N_\mathrm{perRev}(0)\: 
\frac{n_\mathrm{max}(\epsilon)}{n_\mathrm{mean}(\epsilon)} = 
N_\mathrm{perRev}(0)\: \sqrt{\frac{1 + \epsilon}{1 - 
\epsilon}}\:\frac{1}{1 - \epsilon} \ep
\end{equation}
In legends of graphics the quantity $\mathrm{nPerRev}$ always means 
$N_\mathrm{perRev}(0)$ and the number of integration steps per revolution 
that actually is being used depends on $\epsilon$ and is given as 
$N_\mathrm{perRev}(\epsilon)$ from \eq{epsnorm}.
For an estimate of the stability limit it is plausible to take
$\mathrm{nPerRev}$ as the `steps per period' of the surrogate equation
\eq{linear test} and to infer from the analysis there the 
limits $\pi$ for DALF and $\frac{3}{2} \pi$ for ADALF.

\subsection{The main difference between DALF and ADALF}
\label{difference}
\begin{figure}
\centering
\mbox{
\includegraphics[width=120mm]{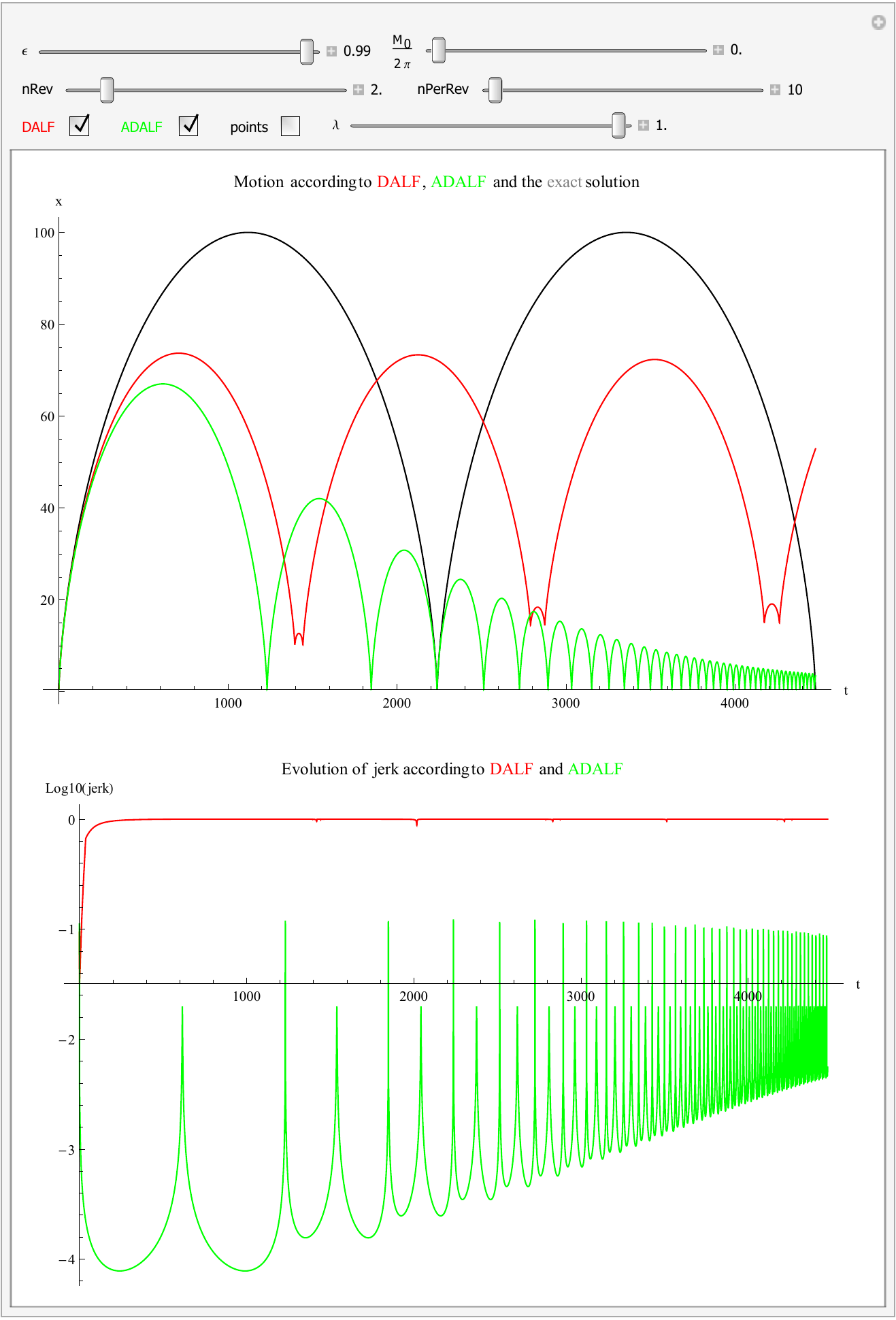}
}
\caption{Position and jerk as functions of time}
\label{fig:jerk1}
\end{figure}

\begin{figure}
\centering
\mbox{
\includegraphics[width=120mm]{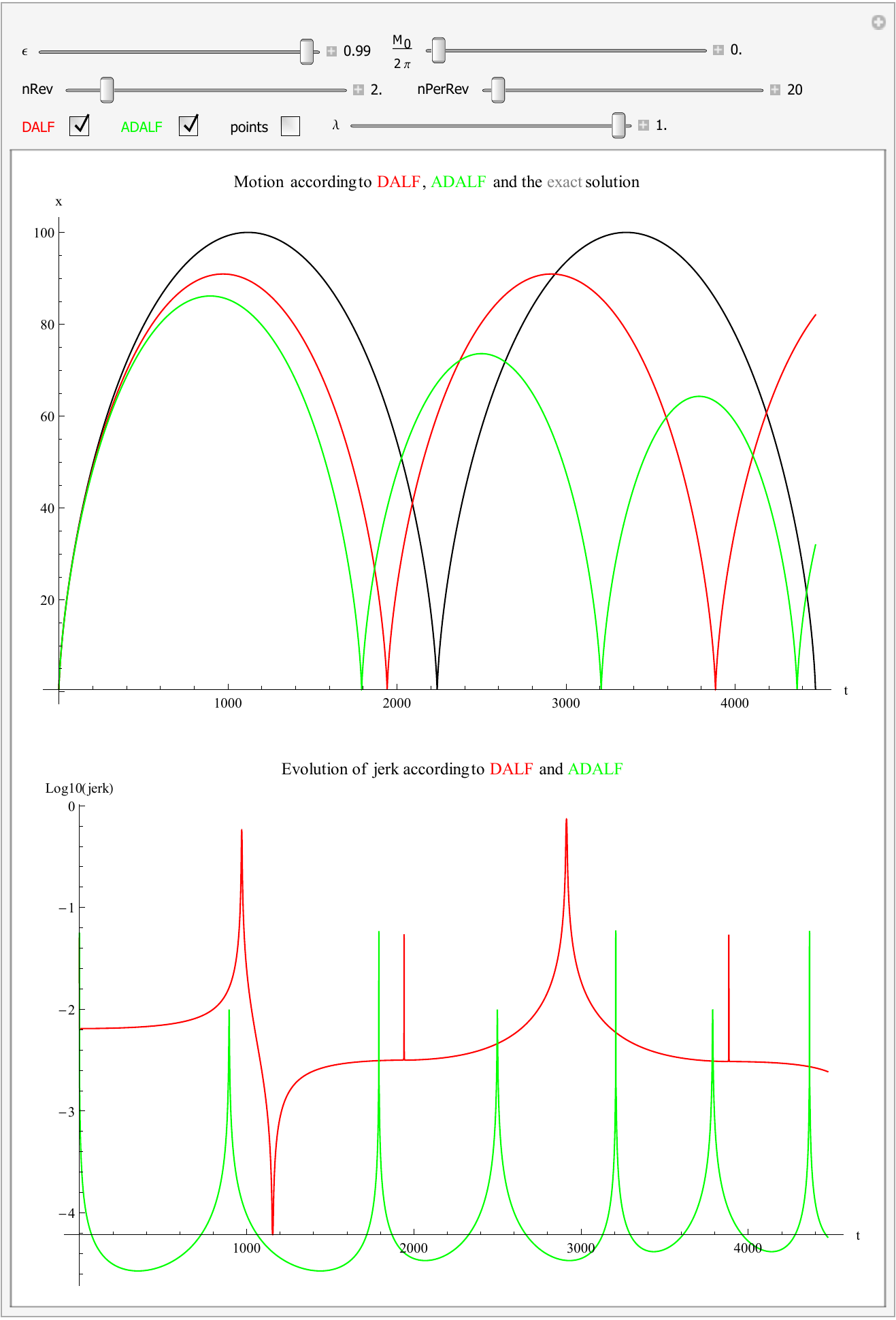}
}
\caption{Position and jerk as functions of time for half of the original 
time step}
\label{fig:jerk2}
\end{figure}

\begin{figure}
\centering
\mbox{
\includegraphics[width=120mm]{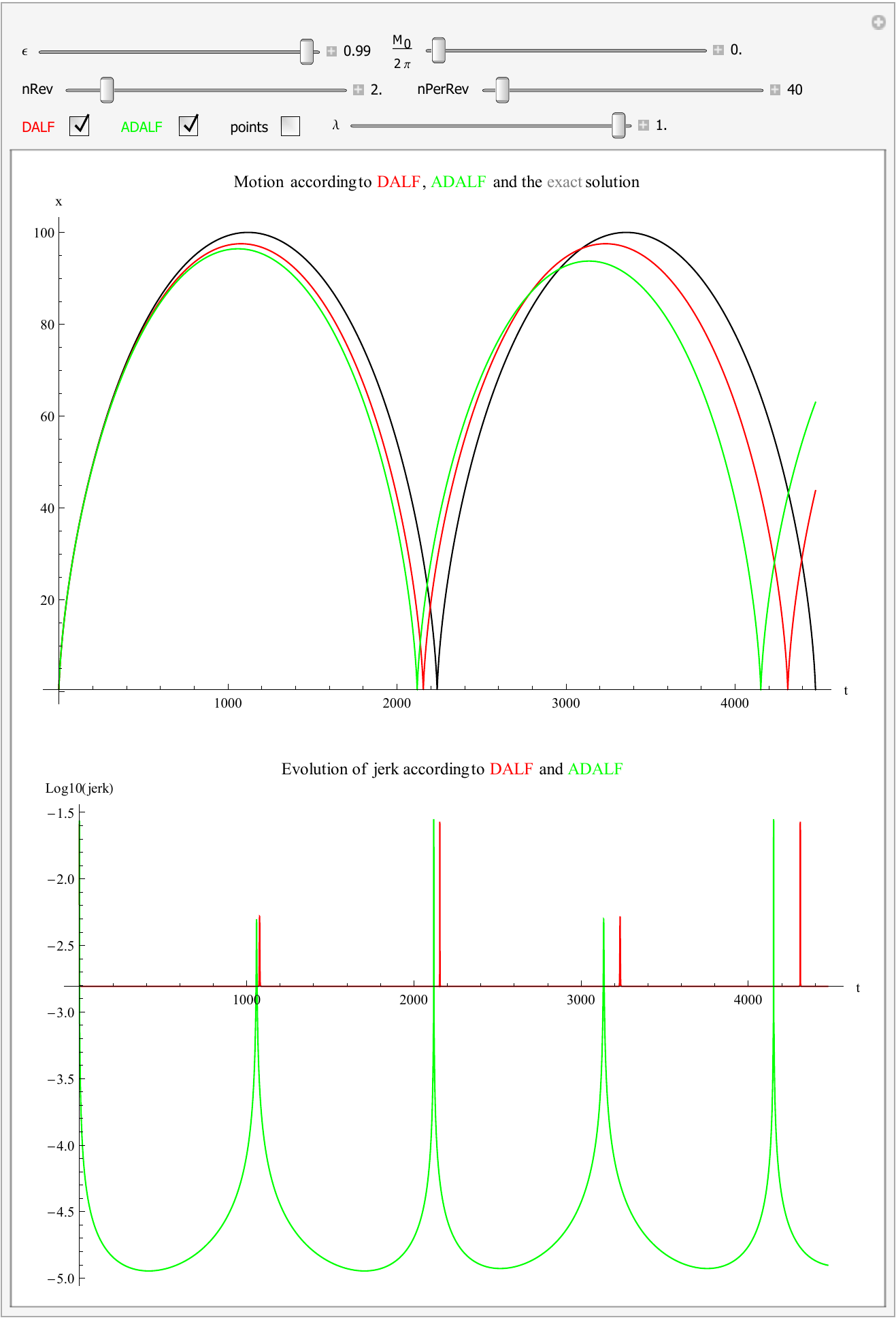}
}
\caption{Position and jerk as functions of time for a quarter of the 
original time step}
\label{fig:jerk3}
\end{figure}

We want to see how DALF's stability deficit manifests 
itself in a typical situation and how ADALF remedies the 
deficit.  
The figures \ref{fig:jerk1}, \ref{fig:jerk2}, and \ref{fig:jerk3} show 
such a situation.
All these figures refer to a trajectory with the extreme eccentricity 
$\epsilon = 0.99$. 
For the first of these figures the time step is by a factor $10 : \pi = 
3.2$  below the stability limit of DALF, so that we are sure that the 
solution will never grow exponentially. 
A certain amount  of detail rendition should not be expected unless the 
time step is by a factor of about $10$ below the stability limit.
In Figure ~\ref{fig:jerk1} we see a dramatic deviation of the DALF 
trajectory from expected behavior. 
Long before having reached its correct turning point $x \approx 0.5$ the
particle stops, makes a feeble leap and then a strong bolt which 
releases all the energy that the particle had prior to the first stop.  
How can the particle be forced to bounce back against the direction of 
the force field?
This is caused by a phenomenon in the working of the DALF integrator that
one may describe as an excitation of an internal degree of freedom.
If this excitation is very large it lets the trajectory attain a zigzag 
character. But long before a zigzag becomes visible in a typical graphical 
representation of a trajectory unexpected things like our bouncing back 
phenomenon may happen.
It is useful to have a quantity which indicates the `degree of internal 
excitation'. Such a quantity can be defined as follows: 
In the algorithm \eq{dalf step} the terms $a := F(t, \psi)$ and $b := 
\phi$ are intended to be approximately equal so that their small 
difference, in 
the third and the sixth step of the algorithm causes a fine tuning of the 
direction of motion in accordance with the vector field $F$. 
To detect deviations from this normal behavior we compute the quantity 
\begin{equation}
\Lb{kappa}
 \kappa(a,b) := \frac{\norm{a - b}}{\norm{a} + \norm{b} +\; \text{tiny}}
\end{equation}
which takes values in $[0,1]$, where values larger than, say, $0.1$ 
indicate significant deviation from $a \approx b$. 
In one DALF step we compute two $\kappa$-values, and their mean 
% define eliminated 2013-09-18
value gives the quantity we wish to define. I refer to this particular 
quantity as \emph{jerk}. 
For ADALF the same definition makes sense and will be used for comparison.
It is clear that the averaging step in ADALF smooths away any zigzag and 
keeps jerk down.
The jerk values for DALF- and for ADALF-trajectories are shown in the 
lower part of the Figures under consideration. 
The first of these figures clearly shows that for the DALF trajectory 
the integrator works fully in the `jerky mode' and thus is no longer under 
the only control of the field $F$ but also reacts to the accumulated 
value of the now virtually autonomous quantity $\phi$. 
The ADALF-trajectory shows considerable jerk only near the return points. 
It absorbs kinetic energy to the effect that the trajectory deviates 
strongly from the exact one, in a predictable manner, though.
The next figures reduce the time step by a factor $2$ each and show no 
longer the unphysical bouncing-back and the trajectories and oscillation 
frequency and amplitude come closer to those of the exact solution. 
The jerk graphs demonstrate that the jerk level goes down systematically.
In the jerk graph of DALF of Figure~\ref{fig:jerk2} there looms a plateau 
which is fully developed in Figure~\ref{fig:jerk3}. Notice that the 
automatically drawn $t$-axis follows just this plateau and reduces its 
visibility.

One should recall that DALF is reversible so that the trajectory, even if 
it has dissolved into a felt of interwoven zigzag lines, can be 
followed back to the stage where it was an innocent smooth curve.
This suggests to consider the `excited state of a leapfrog trajectory' an 
interesting object to study, possibly as a model for a quantum-mechanical 
wave function. 
See also Figure~\ref{fig:ALF relaxation}, case $\lambda = 1$, where the 
% --> , 2013-09-18
excited trajectory mutates back to a smooth one in a periodic manner.  

\subsection{Order}
\label{order}
\begin{figure}
\centering
\mbox{
\includegraphics[width=120mm]{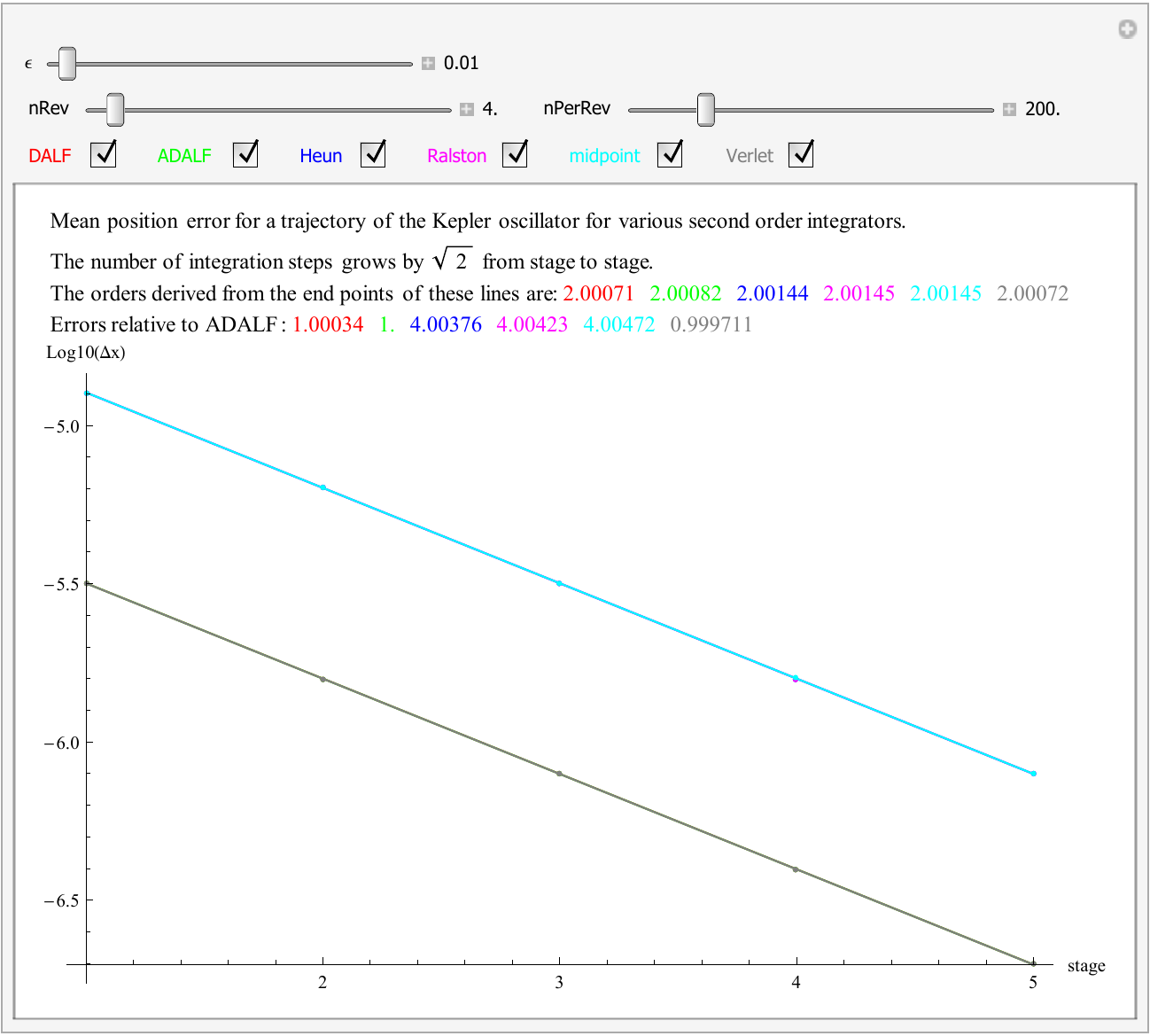}
}
\caption{Order and accuracy of integrators for $\epsilon = 0.01$}
\label{fig:order1}
\end{figure}

\begin{figure}
\centering
\mbox{
\includegraphics[width=120mm]{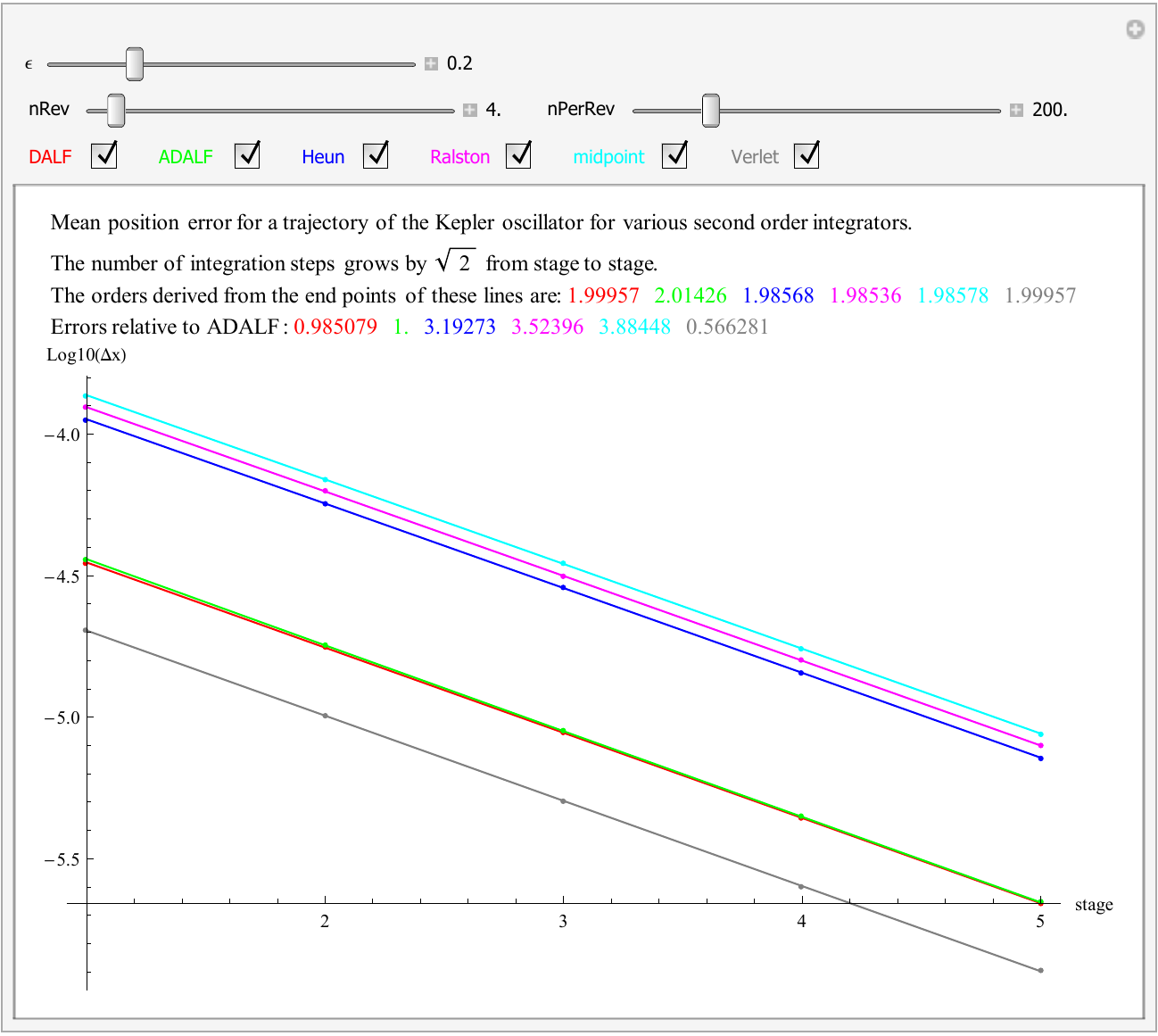}
}
\caption{Order and accuracy of integrators for $\epsilon = 0.2$}
\label{fig:order2}
\end{figure}

\begin{figure}
\centering
\mbox{
\includegraphics[width=140mm]{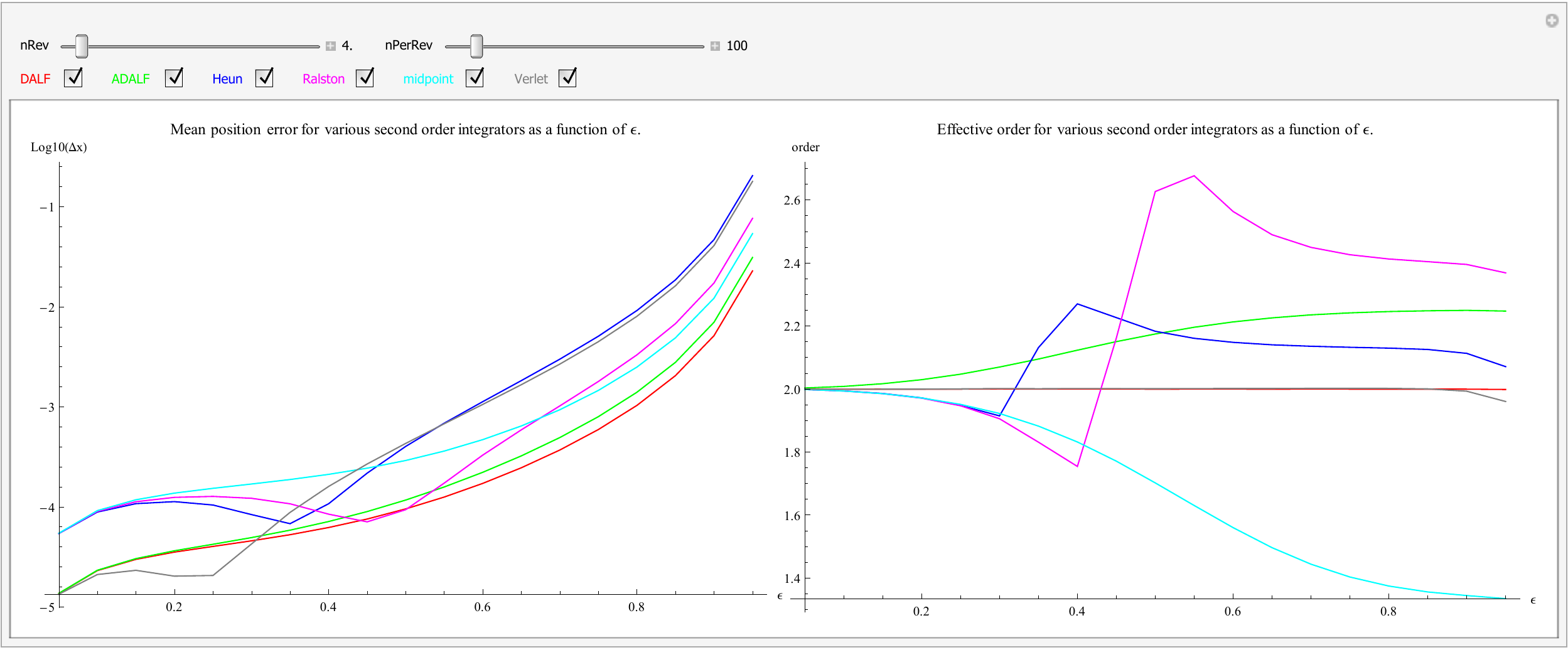}
}
\caption{Accuracy and order of integrators over an $\epsilon$-range 
from $0.05$ to $0.95$ }
\label{fig:order3}
\end{figure}

The order of an integrator is best demonstrated by analyzing how the 
mean error along a trajectory over a given time span depends on the number 
of integration steps. Here, we vary the integration time step over two 
octaves in 4 stages. 
For Figure~\ref{fig:order1} we consider a situation that is 
computationally not demanding: the eccentricity is small ( $\epsilon = 
0.01$) and and the time step is many times below the stability limit for 
% 0.1 replaced by 0.01 2013-09-18
all integrators under consideration. These methods are: DALF, ADALF, the 
three Runge-Kutta methods defined in Subsection~\ref{stability}, and the 
St\o rmer-Verlet integrator which in Section~\ref{Kepler} was referred to 
as direct midpoint integrator.
The result from this figure is that all these integrators are of order $2$ 
and that the leapfrog integrators (to which also St\o rmer-Verlet belongs) 
are four times more accurate as the Runge-Kutta integrators.
For Fig~\ref{fig:order2} the eccentricity is increased by a 
factor of $20$ ($\epsilon = 0.2$) 
% doubled --> increased by a factor of 20 2013-09-18
and we see some differentiation in the data of the methods.
Fig~\ref{fig:order3} shows the mean error and the order for the range $ 
\epsilon = 0.5(0.5)0.95$. 
Here we find variations of the curves from integrator to integrator, 
especially among the Runge-Kutta integrators, which may come as a surprise 
if one recalls that the stability region is the same for the Runge-Kutta 
methods under consideration.
One more observation: The mean error grows with $\epsilon$ although with 
\eq{epsnorm} measures were taken to equalize the accuracy of the result 
with respect to $\epsilon$.

\clearpage

\subsection{Numerical interaction picture}
\label{nip}
\begin{figure}
\centering
\mbox{
\includegraphics[width=140mm]{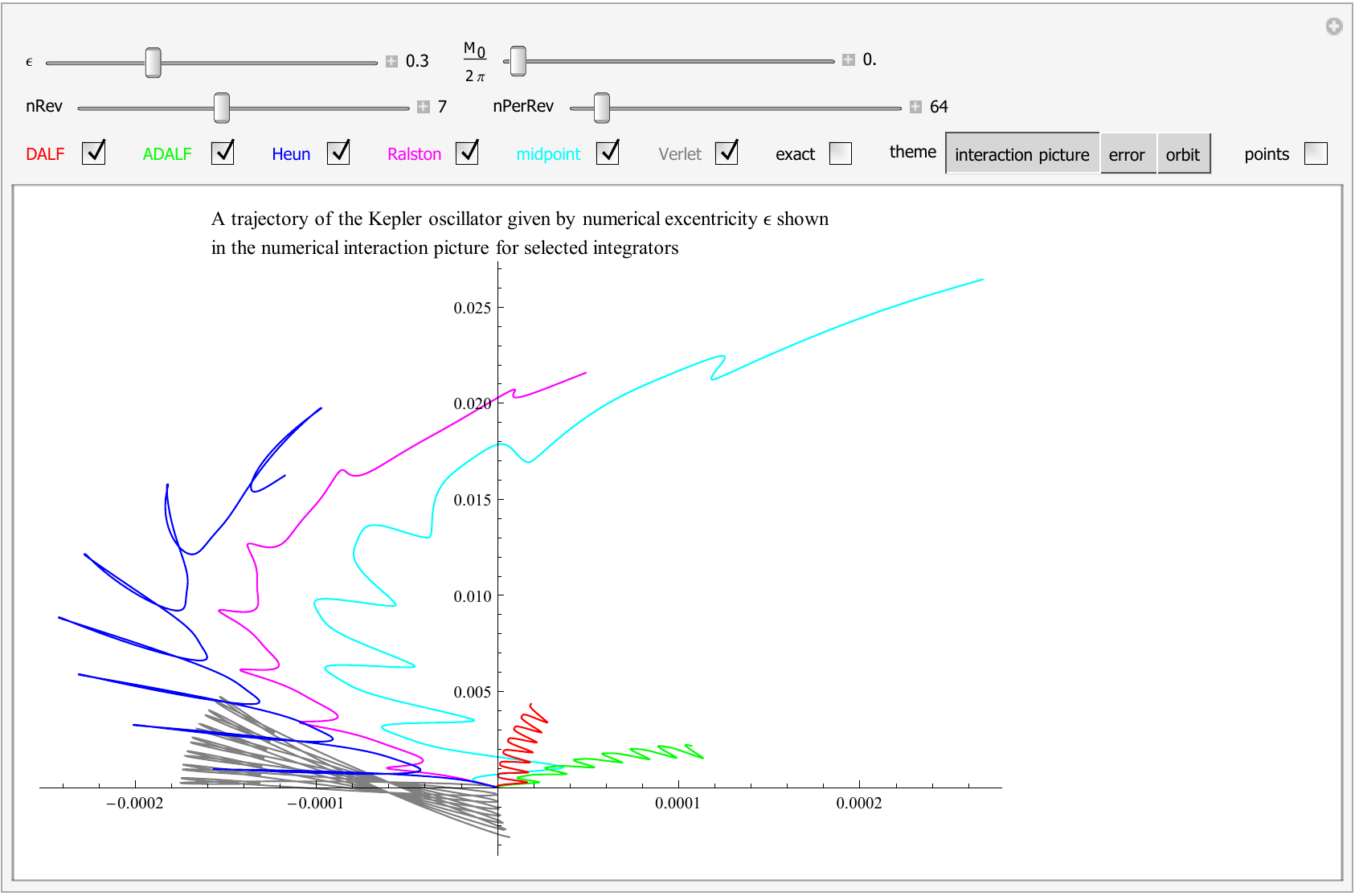}
}
\caption{Interaction picture of Runge-Kutta and leapfrog trajectories}
\label{fig:nip1}
\end{figure}

The trajectories from the numerical interaction picture as introduced in 
Section~\ref{Kepler} show in Figure~\ref{fig:nip1}, in agreement with the 
results of the previous subsection, that the Runge-Kutta methods are 
considerably less accurate than the leapfrog methods. 
Recall that short trajectories indicate high accuracy of the integrator.

\begin{figure}
\centering
\mbox{
\includegraphics[width=140mm]{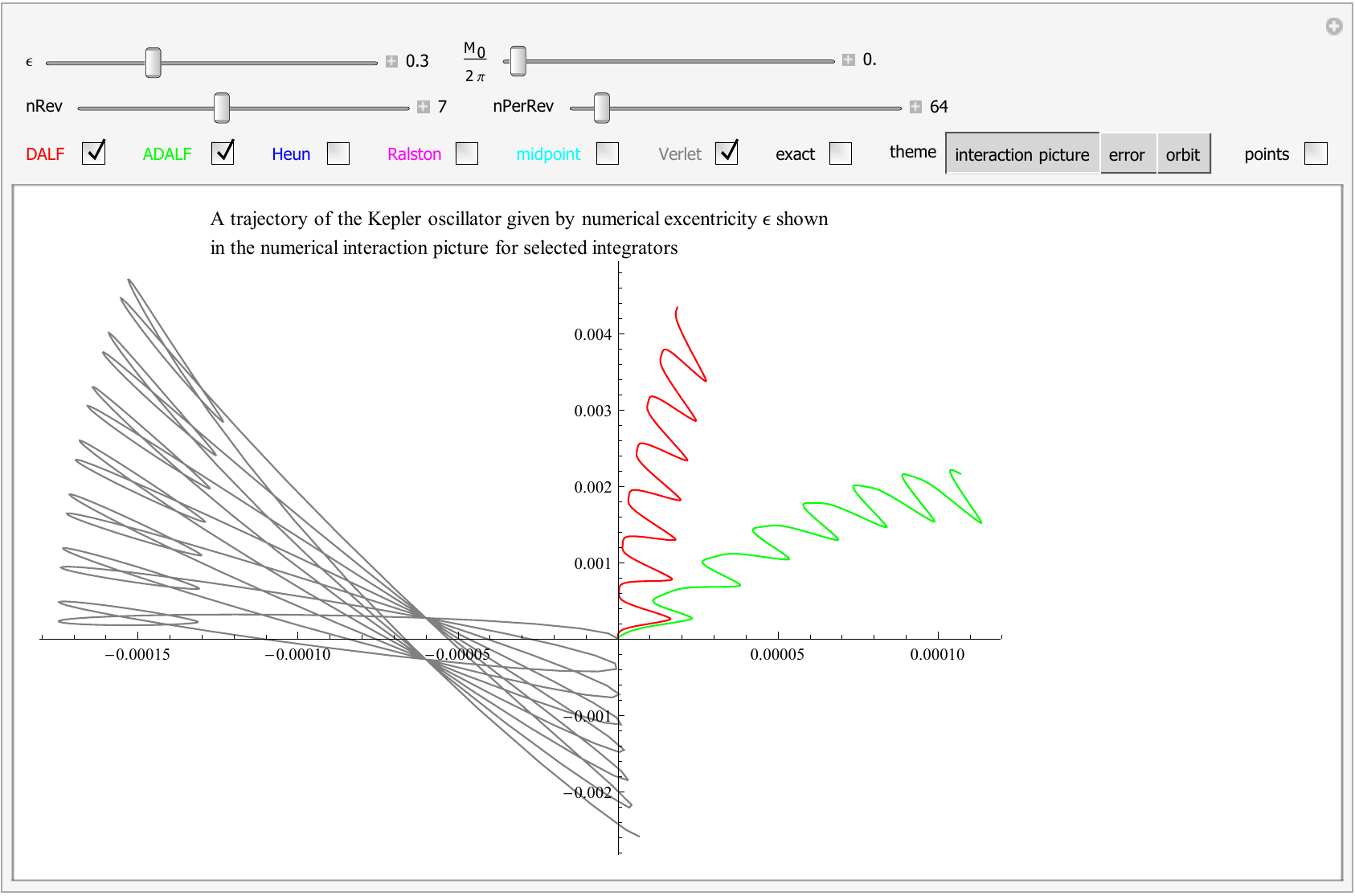}
}
\caption{Interaction picture of leapfrog trajectories}
\label{fig:nip2}
\end{figure}

\begin{figure}
\centering
\mbox{
\includegraphics[width=140mm]{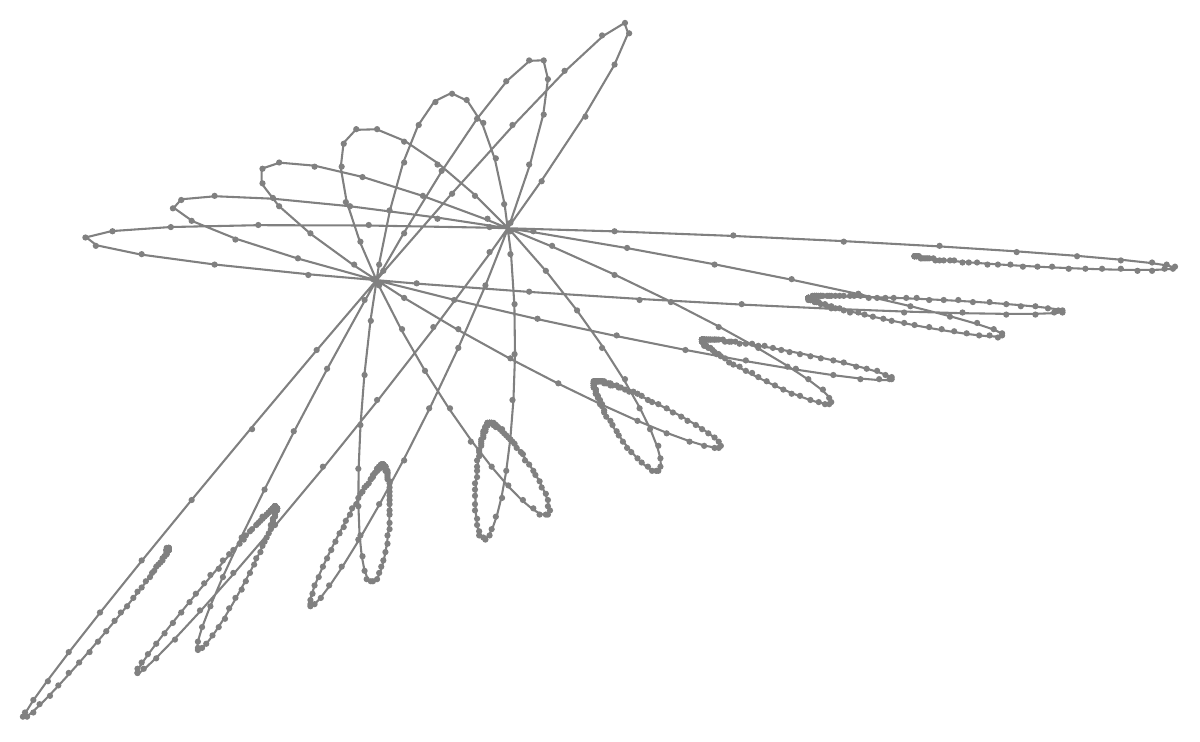}
}
\caption{Interaction picture of a St\o rmer-Verlet trajectory}
\label{fig:nip3}
\end{figure}

Figure~\ref{fig:nip2} shows only the leapfrog methods and thus allows a 
more detailed rendition of their trajectories. 
Here parameters are selected such that in the gray St\o rmer-Verlet 
trajectory a strange `two focal points' feature appears, the nature of 
which is not understood so far.
These focal points are shown again with slightly different parameters in 
Figure~\ref{fig:nip3} in a format that better allows to appreciate the 
interesting geometry and the aesthetic qualities of the trajectory.  

\subsection{Automatic step control}
\label{auto}
\begin{figure}
\centering
\mbox{
\includegraphics[width=140mm]{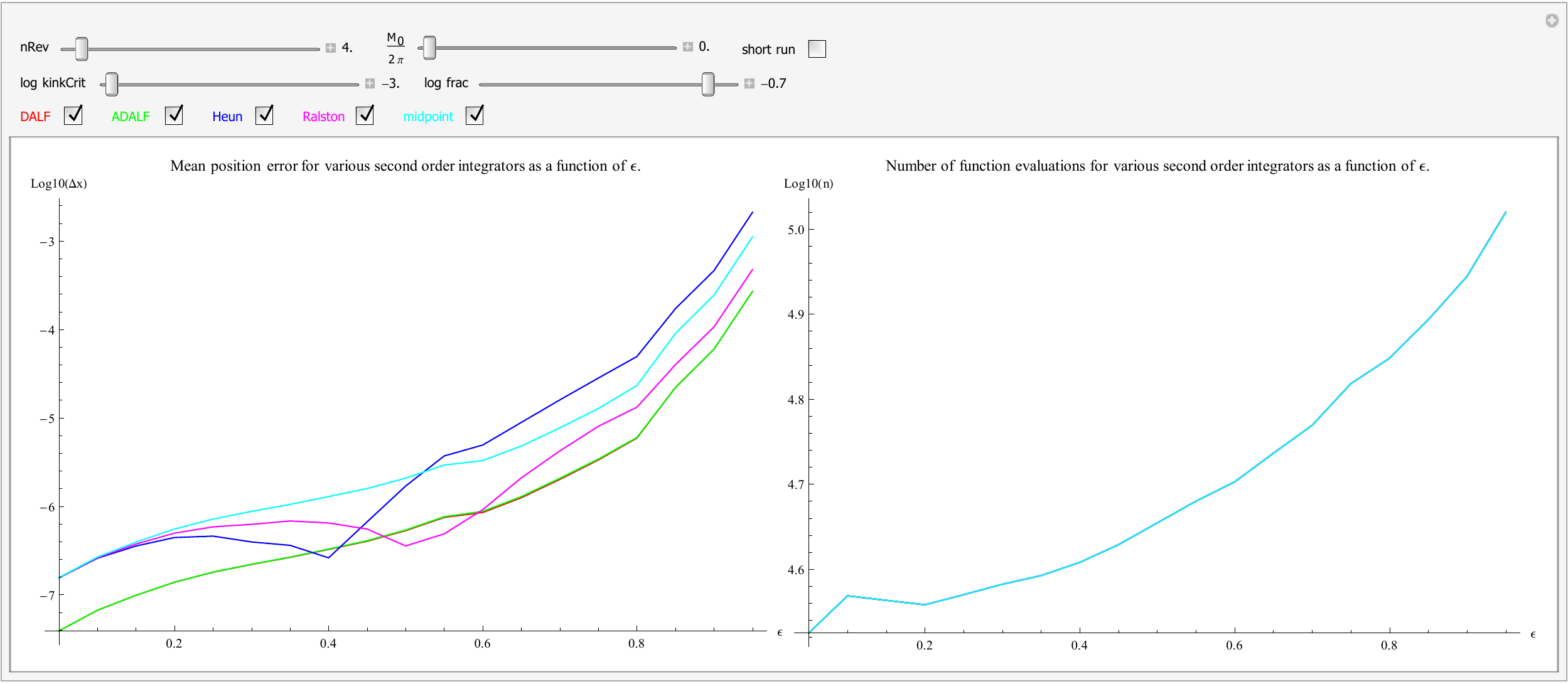}
}
\caption{Run of all integrators, except of St\o rmer-Verlet, under 
automated step control for $\epsilon 
= 0.05(0.05)0.95$}
\label{fig:auto1}
\end{figure}

One of the motivations for inventing the asynchronous version of the 
leapfrog integrator was to facilitate automatic time step control.
This subsection will show that this works in a simple and 
quite effective manner. 
This method seems to have been first demonstrated in \cite{stiff}. 
It assumes that the state to be propagated is of the form $(\psi,\phi)$, 
where $\phi$ was initialized by \eq{initial phi}. 
All integrators which as the first operation in an evolution step do such 
an evaluation of $F$ for values $t$ and $\phi$ valid for the beginning of 
the step can be transformed to this assumed form: 
We leave out this first operation, replace their then no longer available 
result by the available quantity $\phi$ and add at the end of the step an 
operation $\phi = F(t,\psi(t))$ to have $\phi$ available for the step to 
follow.
As already indicated in Subsection~\ref{stability} this can be done for 
all second order explicit Runge-Kutta methods so that we can use those for 
comparison without any change in the driver logic of the program that 
created Figure~\ref{fig:auto1}. 
For the St\o rmer-Verlet method, $\phi$ is not an auxiliary quantity that 
can be given a value by evaluation of $F$. So the method does not 
apply verbally. So it is left out as a comparison method here.

So, how time step control works?
After having done an evolution step we 
know an initial value $\phi_\text{i}$ and a final value $\phi_\text{f}$ of 
$\phi$. 
Using \eq{kappa} we compute $\kappa(\phi_\text{i},\phi_\text{f})$
(this is a simple computational step, which does not require 
an evaluation of $F$).
If this dimensionless quantity exceeds some agreed critical level $a_1$ we 
conclude that this step was not acceptable. 
As a consequence we re-initialize $\phi$ and re-do the step with a 
shorter time step. Let us call $f_1$ the factor with which we multiply the 
old time step. 
If $\kappa$ is smaller than some agreed comfortable level $a_2$, we 
increase the time step to be used for the next evolution step by 
multiplication with some some agreed factor $f_2$.
In the program of Figure~\ref{fig:auto1} there are input quantities 
kinkCrit = 0.001 and frac = 0.2. 
These determine the quantities just mentioned as follows: 
\[ a_1 = \text{kinkCrit} \ec a_2 = 0.5 a_1 \ec f_1 = 1 - \text{frac} \ec 
f_2 = 1 + \text{frac} \ep \]
We see that the results are by a factor of about 200 more accurate than in 
the fixed step integration Figure~\ref{fig:order3}. However we need 25 
times more function evaluations. 
For a second order integrator this translates into a factor 625 for 
accuracy which has to be compared to the factor 200 that we actually got. 
So the method is by a factor 3 less efficient than one would expect for an 
optimally selected constant time step.
Of course, all integrators need nearly the same number of function 
evaluations 
since curvature of the exact trajectory determines, where a 
re-initialization is needed.

Over the $\epsilon$-range in Figure~\ref{fig:order3} the time 
step varies by a factor of about 125 and 
only of about 3 in Figure~\ref{fig:auto1}.

The error curves in Figure~\ref{fig:order3} and Figure~\ref{fig:auto1} 
show surprising similarity. 
Especially apparent are the minima of the curves of  Ralston's and Heun's 
methods for $\epsilon \approx 0.4$.

For serious applications in which both accuracy and computation time 
matter, modified values of $a_1, a_2, f_1, f_2$ or even variations of the 
algorithm may be useful. 
My experience with large applications is restricted to constant time step.

\clearpage

\section*{Acknowledgment}
I am grateful to Domenico Castrigiano for many discussions on the relation 
of discrete mathematics to classical analysis and to Ernst Hairer , 
Blair Perot, and David Seal for valuable comments. H.E. Lehtihet has 
provided instructive test examples and provided ideas for
interpreting the obtained results. 
Sajad Jafari shared his very interesting and computationally demanding 
equations dealing with `romantic relationships' prior to publication. 
Finally I am grateful to J.M. Sanz-Serna for having made me aware of the
Nordsieck technique and Butcher's systematic of Linear General Methods.

last modification 2016-04-14

\end{document}